\documentclass[a4paper]{article}

\usepackage{latexsym}
\usepackage{amsbsy}
\usepackage{amsmath}
\usepackage{amssymb}
\usepackage{graphics}
\usepackage{epsfig}

\def\ZZ{\mathbb{Z}}

\def\i{\mathrm{int}}

\newtheorem{thm}{Theorem}[section]
\newtheorem{lem}[thm]{Lemma}
\newtheorem{prop}[thm]{Proposition}

\newenvironment{defin}
      {\par\noindent 
                \textbf{Definition.}}
             {}

\newenvironment{proof}%
      {\par\noindent%
            \textbf{Proof.}}%
           { ~\hfill$\Box$\linebreak}

\newcommand{\po}{ (1\> - \>4) }
\newcommand{\pt}{ (2\> - \>3) }
\newcommand{\pth}{ (3\> - \>2) }
\newcommand{\pf}{ (4\> - \>1) }
\newcommand{\p}{ \> - \> }

\title{\large \bfseries SIMPLIFYING TRIANGULATIONS OF 
                                \boldmath $S^3$}

\author{\normalsize ALEKSANDAR MIJATOVI\'{C}}

\date{}

\renewenvironment{abstract}
            {\begin{quotation}\noindent\small 
                              \textsc{Abstract}.\hspace{0.5mm}}
            {\end{quotation}}

\begin{document}

\maketitle

\begin{abstract}
In this paper we describe a procedure to simplify any given
triangulation of $S^3$ using Pachner moves. We obtain an explicit 
exponential-type bound
on the number of Pachner moves needed for this process. This leads to a new 
recognition algorithm for the 3-sphere.
\end{abstract}

\begin{center}
\section{\normalsize \scshape INTRODUCTION}
\label{sec:intro}
\end{center}

It has been known for some time that any triangulation of a 
closed PL 
$n$-manifold
can be transformed into any other triangulation of the same manifold
by a finite sequence of moves~\cite{pachner}. We can describe the moves 
as follows.\\

\begin{figure}[!hbt]
  \begin{center}
%     \vspace{2cm}
%     \Huge{(2-3) (1-4)}
    \epsfig{file=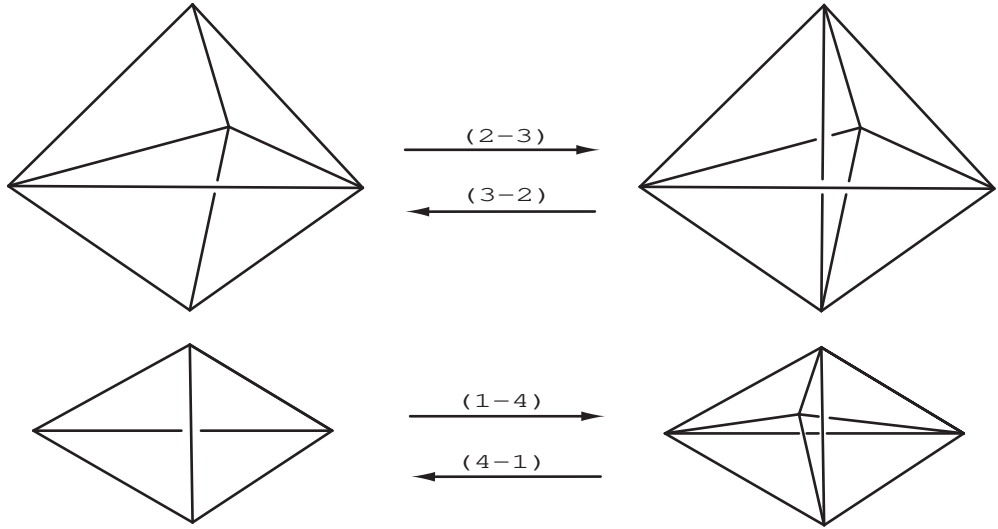}
        \caption{{\small Three dimensional Pachner moves.}}
        \label{fig:3pm}
   \end{center}
\end{figure}

\begin{defin}
Let 
$T$
be a triangulation of an
\mbox{$n$-manifold}
$M$.
Suppose
$D$
is a combinatorial 
\mbox{$n$-disc}
which is a subcomplex both of
$T$
and of the boundary of a standard 
\mbox{$(n+1)$-simplex}
$\Delta^{n+1}.$
A
\textit{Pachner move}
consists of changing
$T$
by removing the subcomplex
$D$
and inserting
$\partial\Delta^{n+1}-\i(D)$
(for 
$n$
equals 3, see figure~\ref{fig:3pm}).
\end{defin}\\

It is an immediate consequence of the definition that there are precisely
$(n+1)$
possible Pachner moves in dimension
$n$.
We can now state
Pachner's result~\cite{pachner} in the following way.\\

\begin{thm}[Pachner]
\label{thm:pach}
Closed PL 
$n$-manifolds
$M$
and
$N$,
triangulated by
$T$
and 
$K$
respectively,
are piecewise linearly homeomorphic
if and only if
there exists a finite sequence of Pachner moves
and simplicial isomorphisms
taking the triangulation 
$T$
into the triangulation 
$K$.
\end{thm}

In dimension 3
we have 
four moves from figure~\ref{fig:3pm} at our disposal.
Using them, we can describe the main theorem of this paper.

\begin{thm}
\label{thm:main}
Let 
$T$
be a triangulation of a 3-sphere and let 
$t$
be the number of tetrahedra in it. Then we can simplify
the triangulation 
$T$
to the canonical triangulation 
of
$S^3$,
by making less than
$a\>t^22^{bt^2}$
Pachner moves, where the constant
$a$ 
is bounded above by 
$6\cdot10^{6}$
and the constant
$b$
is smaller than
$5\cdot10^4$.
\end{thm}

The triangulation 
$T$
in this theorem can be non-combinatorial
(i.e. simplices are not uniquely determined 
by their vertices), as is the case with
the canonical triangulation of
$S^3$,		
consisting of two standard 3-simplices glued together via an
identity on their boundaries. We should mention here that 
the original Pachner's proof of the theorem~\ref{thm:pach}
works for combinatorial triangulations only. However, at least
in dimension 3, this does not matter because the second 
derived subdivision of any (possibly non-combinatorial)
triangulation is always combinatorial and can be obtained from the
original triangulation by a finite sequence of Pachner moves. 

A possible effect Pachner's result could have on the
theory of 
\mbox{$3$-manifolds}
is discussed by the next proposition.

\begin{prop}
\label{prop:computable}
Let
$T$
and
$K$
be two triangulations of the same
closed PL
\mbox{$3$-manifold}
$M$.
The
existence of a computable function, depending only on the number
of
\mbox{$3$-simplices}
in 
$T$
and
$K$,
bounding the number of Pachner moves required to transform 
$T$
into
$K$,
is equivalent to an
algorithmic solution of the recognition problem for
$M$
among all
\mbox{$3$-manifolds}.
\end{prop}

\begin{proof}
Assume first that 
$f(t,k)$
is a computable function as described in the proposition. 
Suppose that 
$T$
is a triangulation of 
$M$
with 
$t$
3-simplices. 
Let
$K$
be a triangulation of some closed
\mbox{$3$-manifold}
$N$
containing 
$k$
3-simplices.
Do all possible sequences of Pachner moves 
on the triangulation 
$T$
of length
at most
$f(t,k)$,
and check each time if the result is isomorphic to
$K$.
This gives an algorithm to determine whether 
$M$
and
$N$
are PL homeomorphic.

Conversely suppose that we have an algorithm to 
recognize 
$M$ among all 3-manifolds.
Now we need a complete (finite) list 
of all triangulations of all 
\mbox{$3$-manifolds}
with 
a fixed number of 
\mbox{$3$-simplices}.
In dimension three, such a list can be built algorithmically
because there is an easy way of recognizing the 2-sphere 
(the Euler characteristic suffices) as a link of a vertex.

We can now create 
all triangulations of 
$M$
with the specific number of 3-simplices by running the 
recognition algorithm for 
$M$
(which exists by assumption) 
on the list of all 3-manifold triangulations
with the specified number of 3-simplices.

An algorithm, making all possible Pachner moves on a triangulation
of our
\mbox{$3$-manifold}
$M$
with
$t$
3-simplices will
after a finite number of steps 
(by the theorem~\ref{thm:pach}) necessarily produce
a given triangulation of
$M$
containing
$k$
3-simplices.
Since we can list all triangulations of
$M$
with 
$t$
(respectively 
$k$)
3-simplices, this gives an algorithm to calculate the value
of
the function
$f(t,k)$
as required.
\end{proof}

At present there is no known algorithm 
to decide
whether a given simplicial complex  
is an
$n$-sphere,
for 
$n\geq 4$.
This means that the proof of one 
of the implications in proposition~\ref{prop:computable}
breaks down in dimensions
five and above since there is no way of 
building a list of all triangulations of all manifolds
with a fixed number of top dimensional simplices, in 
these dimensions.

The proof of the
converse implication in proposition~\ref{prop:computable},
showing that a computable bound implies a recognition
algorithm for a given
$n$-manifold,
remains valid in any dimension.
Furthermore, if such a computable bound 
existed 
for all 
$n$-manifolds,
and was independent of the underlying
$n$-manifold,
then it would give an algorithm to determine
whether any two 
$n$-manifolds are homeomorphic.
But using the fact (proved by A.A. Markov) that there is no such
algorithm for 
$n\geq 4$,
we can conclude that such a computable function 
does not exist in dimensions four and above. 

It is interesting to note, that 
for any 
$n$-manifold 
$M$ 
Pachner's theorem implies the
existence of a function, depending only on the number of
$n$-simplices
in
$T$
and
$K$,
and bounding the number of Pachner moves necessary for the 
whole transformation. This is because we can generate 
(in principle, at least) all possible triangulations of
our 
$n$-manifold
$M$
with fixed numbers of 
$n$-simplices. Then, using theorem~\ref{thm:pach}, 
a finite sequence of Pachner moves connecting any two of
them, can be found. Taking the maximum length over this
finite family of sequences gives us the bound. 
Therefore, computability of the function in 
proposition~\ref{prop:computable} is an assumption that can
not be omitted.

The upper bound in theorem~\ref{thm:main} is computable.
It therefore yields
a new recognition 
algorithm for the 
3-sphere.

\begin{center}
\section{\normalsize \scshape NORMAL SURFACES}
\label{sec:normal}
\end{center}

One of the essential ingredients of the proof of 
theorem~\ref{thm:main} is the theory of normal 
and almost normal surfaces. In
this section we shall describe some of its basic features. We will
then go on to discuss the Rubinstein-Thompson 
algorithm~\cite{thompson} 
for recognizing the 3-sphere which provides the setting for the
proof of theorem~\ref{thm:main}. 
After it, we'll mention some of the consequences
of normal surface theory which will prove to be useful later.
At the end of this section we shall prove the isotoping
lemma that will later give us a way of simplifying triangulations
of the 3-sphere.
Let's start with some definitions.

A
\textit{normal triangle}
(respectively \textit{quadrilateral})
in a 3-simplex 
$\Delta^3$
is a properly embedded disc
$D$,
such that its boundary
$\partial D$
intersects precisely three (respectively four)
edges transversely in a single point and is disjoint from
the remaining 1-simplices and vertices of 
$\Delta^3$.
A 
\textit{normal disc}
is a normal triangle or quadrilateral.

\begin{figure}[!hbt]
  \begin{center}
%     \vspace{2cm}
%     \Huge{NORMAL DISCS}
    \epsfig{file=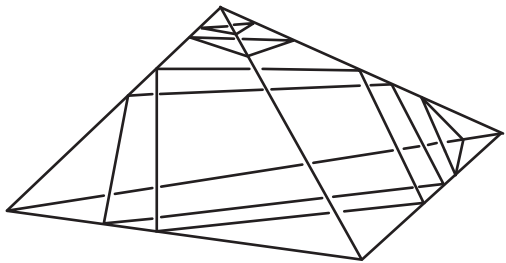}
        \caption{{\small Three types of normal discs.}}
        \label{fig:normal}
   \end{center}
\end{figure}

There are four possible types of normal triangles, because
each triangle is parallel to one of the faces of
$\Delta^3$. Normal quadrilaterals will always separate the 
vertices of the tetrahedron in pairs. It is therefore clear,
that we can only have three possible quadrilateral types. Together,
there are 7 distinct normal disc types in a tetrahedron.

Let 
$M$
be a 3-manifold with a triangulation 
$T$.
A properly embedded surface 
$F$
in 
$M$
is in 
\textit{normal form}
with respect to the triangulation
$T$,
if it intersects each tetrahedron of
$T$
in a finite (possibly empty) collection of disjoint normal 
discs.

Since normal surfaces are always embedded, at most one 
of the quadrilateral types can occur in each 3-simplex.

Suppose
$F$ 
is a normal surface in 
$M$
with respect to 
$T$.
Then
$F$
corresponds to a vector
$\underline{x}=(x_1, \ldots ,x_{7t})$
with 
$7t$
coordinates, where 
$t$
denotes the number of 3-simplices in the triangulation
$T$.
The index set
$\{1, \ldots ,7t\}$
corresponds to all possible disc types in 
$T$ 
(there is 7 of them for each tetrahedron). 
The coordinate
$x_i$
is simply the number of copies of 
\mbox{$i$--th}
disc type in our surface
$F$.

Each 2-simplex in 
$T$
contains three types of normal arcs (coming from 
normal discs), one cutting off each vertex of the triangle.
If it is a face of two 3-simplices in
$T$,
then it gives rise to three matching (linear) equations, one
corresponding to each normal arc type. 
Doing this for every triangle, not in the boundary of
$M$,
we've constructed a linear system 
in
$7t$
variables,
consisting of
at most 
$6t$
equations.

It follows immediately from the construction, that 
the vector 
$\underline{x}$,
coming from the normal surface 
$F$,
gives a solution to the linear system. By
imposing extra conditions to ensure that all
quadrilaterals in a given tetrahedron are of the
same type, we obtain a restricted linear system. 
The conditions we've just added are sometimes
referred to as
\textit{quadrilateral constraints}.
Now there is a one to one correspondence between 
embedded normal surfaces in
$M$
and non-negative integral solutions to the restricted
linear system.

Haken proved that all non-negative integral solutions to
such a system are integer linear combinations of a 
finite set of non-negative integral solutions 
$\underline{x_1}, \ldots ,\underline{x_n}$,
called 
\textit{fundamental solutions}, which can be found
in an algorithmic way. As it turns out, these 
fundamental solutions are characterized by the property
of not having a decomposition as a sum of two 
(non-trivial) non-negative integral solutions to the restricted
linear system.

Since each fundamental solution corresponds to an 
embedded normal surface, we obtain a finite set 
$F_1, \dots ,F_n$
of embedded normal surfaces, called
\textit{fundamental surfaces}.
Any embedded normal surface in
$M$
can thus be written algebraically as a non-negative integer
linear combination of fundamental surfaces. Miraculously, this
algebraic fact carries over to geometry. In other words, we 
can define a geometric addition for any two normal surfaces
$F$
and
$G$
with the property that the sum of the corresponding 
solutions to the restricted linear system, is again a 
solution of the same system. This condition boils down
to the fact that the union of all normal discs in both
$F$
and
$G$
satisfies the quadrilateral constraints.

Assuming that and putting both surfaces in general position
with respect to one another, cutting along the arcs of
intersection in each tetrahedron, and pasting the pieces back
together in the unique way, so that we end up with normal discs
only, yields a well-defined embedded normal surface 
$F+G$.
Its corresponding vector is a sum of the vectors coming 
from 
$F$
and
$G$.
The cut and paste process described above is sometimes called
\textit{regular alteration}.

An isotopy of the ambient manifold, preserving the normal 
structure of a given normal surface is called a
\textit{normal isotopy}. 
We should also note that the geometric 
addition described above is well-defined up to 
a normal isotopy of the summands.

Before we describe the Rubinstein-Thompson algorithm, we need to 
introduce a new concept.

\begin{figure}[!hbt]
 \begin{center}
%  \vspace{2cm}
%  \Huge{ALMOST NORMAL}
    \epsfig{file=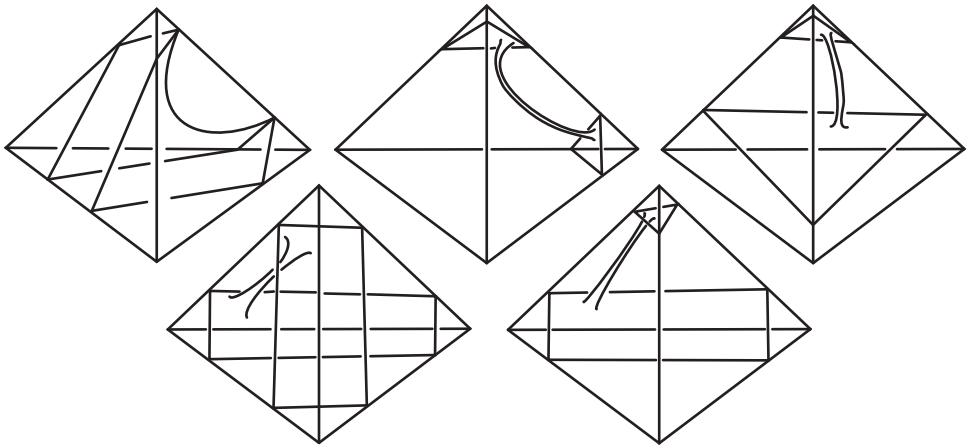}
  \caption{\small Almost normal pieces.}
  \label{fig:almostnormal}
 \end{center}
\end{figure}

\begin{defin}
A properly embedded surface in a 3-manifold
$M$
with a triangulation
$T$
is
\textit{almost normal}
with respect to 
$T$,
if it intersects each tetrahedron of
$T$
in a finite (possibly empty) collection of disjoint normal discs
except in precisely one tetrahedron there is precisely one exceptional
piece from figure~\ref{fig:almostnormal} and possibly some 
normal triangles.
\end{defin}

This exceptional piece is either a disc (the first
possibility on figure~\ref{fig:almostnormal}) whose boundary is 
a normal curve of length eight (i.e. an octagon), or it
is an annulus consisting of two normal disc types with
a tube between them that is parallel to an edge of the
1-skeleton.

Now we can describe the Rubinstein-Thompson algorithm 
which is designed
to determine whether or not a 3-manifold
$M$
with a triangulation
$T$
is a 3-sphere. We can assume that 
$M$
is closed, orientable and that
$H_1(M;\ZZ_2)$
is trivial. All these properties can be checked algorithmically. The
last assumption guarantees that 
$M$
contains no closed non-separating surfaces. 
The algorithm now is in three steps. 
We proceed as follows.

\begin{list}{}{\setlength{\rightmargin0}{\leftmargin}}
\item[\textit{Step} 1.] Find a maximal collection
                        $\Sigma$
                        of disjoint non-parallel normal 2-spheres in
                        $M$.

\item[\textit{Step} 2.] Cut 
                        $M$
                        open along 
                        $\Sigma$. 
                        This splits 
                        $M$
                        into three different types of pieces:
                        \begin{list}{}{\setlength{\rightmargin0}
                                                 {\leftmargin}}
                        \item[\textit{Type} A:]
                                      a 3-ball neighborhood of a vertex of
				      $T$ (every vertex is enclosed in such 
                                      a piece).	
                                      
                        \item[\textit{Type} B:]
                                      a piece with more than one boundary 
                                      component.

                        \item[\textit{Type} C:]
                                      a piece with exactly one boundary
                                      component which is not of type A.
                        \end{list}

\item[\textit{Step} 3.] Search each type C piece for an almost normal 
                        2-sphere with an octagonal component.

\item[\textit{Conclusion}:] $M$
                            is a 3-sphere if and only if every 
                            type C piece contains an 
                            almost normal 2-sphere with an octagonal
                            component.
\end{list}

The bulk of the proof that this indeed is a recognition 
algorithm for the 3-sphere relies on the following two lemmas
from~\cite{thompson}.

\begin{lem}
\label{lem:t1}
A type \textrm{B} piece is a punctured 3-ball.
\end{lem}

\begin{lem}
\label{lem:t2}
A type \textrm{C} piece is a 3-ball if and only if it contains
an ``octagonal'' almost normal 2-sphere.
\end{lem}

By lemma~\ref{lem:t2}, if some type C piece fails to contain 
an ``octagonal'' almost
normal 2-sphere, then it is not a 3-ball and
$M$
is not a 3-sphere. Otherwise, 
$M$
is just a collection of 3-balls and punctured 3-balls glued 
together. Since every 2-sphere is separating, 
$M$
has to be a 3-sphere.

The difficult part of the argument is in the proof of
lemma~\ref{lem:t2}. 
It is here that Thompson simplified Rubinstein's original 
methods to 
prove the existence
of an ``octagonal'' almost normal 
2-sphere in a 3-ball of type C, by 
using Gabai's powerful notion of thin position.
We should also note that the easier converse implication in 
lemma~\ref{lem:t2} follows from lemma~\ref{lem:iso}. 

In order to be able find a maximal collection 
of disjoint non-parallel normal 2-spheres in
$M$
in an algorithmic way, we need the following lemma.

\begin{lem}
\label{lem:maximal}
A maximal collection 
$\Sigma$
of disjoint non-parallel normal 2-spheres in 
$M$,
as in the Rubinstein-Thompson algorithm,
can always be constructed algorithmically.
\end{lem}

A detailed proof of this fact is given
in~\cite{king} (see lemma 3). 
It uses the notion of 
normal surfaces in handle decomposition
that are (loosely speaking) dual to a given triangulation
of our manifold. But since this is only
a technical point that we are clarifying
here, we are not going
to pursue this dual theory any further.  

It will however 
be important for us to be able to bound the complexity
of all of these 2-spheres in
$\Sigma$. 
We shall therefore give a brief description 
of this algorithm. 

Additivity of Euler characteristic implies at once that 
if there exists a 
non-trivial normal 2-sphere in our triangulation,
we can also find one 
(which is also non-trivial)
among fundamental surfaces.
Since the family
of fundamental surfaces is accessible in an algorithmic way, we can
take this fundamental 2-sphere to be the first element in
$\Sigma$.

Now we cut along it. The original triangulation
$T$
of 
our manifold
gives rise to a cellular decomposition of the 
complementary components. What we want to do at this stage
is to look for new normal 2-spheres in each piece separately. 
In order to do that, we need to set up a normal surface theory
in this special kind of cellular complex. It is plausible and 
it turns out to be true that only the cells which are not
bounded by parallel normal pieces of the 2-sphere we cut along,
are the ones where the normal structure of our future surfaces 
can be 
varied. 
This is essentially because in the ``parallel'' regions, a normal
surface has no choice but to run parallel to the boundary.

We can now use the same argument as before to search (in an 
algorithmic way) for new fundamental normal 2-spheres in each 
complementary piece.
A vital fact we are relying on here is that normal 
surfaces in the cut open manifold are also normal in 
our original triangulation. 

So in order to find a maximal family 
$\Sigma$
of disjoint non-parallel normal 2-spheres,
we just have to keep repeating this procedure 
until none of the pieces
in the cut open manifold contain another normal 2-sphere
that is not normally parallel to a boundary component. 
Lemma~\ref{lem:sphere} implies that this process has to terminate.

As far as the complexity, i.e. the number of normal pieces,
of elements in 
$\Sigma$
goes, at each stage it is going to be bounded by 
proposition~\ref{prop:hass}. This bound will
of course be relative to the stage we are at. Since
the linear algebra in the proof of proposition~\ref{prop:hass} 
(which can be found in~\cite{hass}) depends only on the
number of normal variables, it is 
the number of non-trivial cells 
in the complementary components of any normal 2-sphere 
we cut along, that needs to be controlled. The proof
of lemma~\ref{lem:sphere} shows that this number is bounded 
linearly by the number of tetrahedra in
$T$.
We will calculate explicit upper bounds later on in this section.

We still need to answer the question of how to search for
``octagonal'' almost normal 2-spheres in type C pieces. 
Modified versions
of standard normal surface theory algorithms suffice for the
search. Again all the technical details can be found
in~\cite{king}. So our goal is to construct 
an algorithmic procedure which will find
an ``octagonal'' almost normal 2-sphere
in each type C piece. These 2-spheres
will exists by lemma~\ref{lem:t2}
if the 3-manifold 
$M$ 
we are looking at is a 3-sphere. We proceed
as follows.

First we fix a tetrahedron
$H$
in the triangulation
$T$
of a 3-manifold
$M$
and then we fix a normal curve 
$c$
of length eight on its boundary (there are three choices
for
$c$).
Now an analogue to the classical normal surface theory can be
set up. The matching conditions will look 
just like before. Quadrilateral
constraints have to be modified however, 
because we want our solutions to
consist of normal triangles and quadrilaterals everywhere except 
in 
$H$,
where we want them to be composed of normal triangles and 
octagonal components with boundaries parallel to
$c$.
The notion of regular alteration can be defined in this
generalized setting and again it gives rise to
the correspondence between integer linear combinations
of the fundamental solutions to the
(generalized) restricted linear
system and the set of all surfaces described above. 
Fundamental surfaces are again the ones corresponding to 
fundamental solutions. We should also note that their 
complexity is bounded by proposition~\ref{prop:hass} since
the linear system they are the solutions of, has less then 
$7t$
variables.

What we really want is to find the 
``octagonal'' almost normal 2-spheres that are contained 
in type C pieces. The unique boundary component of such
a piece is obtained as a normal 2-sphere in some
cellular decomposition of some piece of the manifold 
$M$.
After we cut along this last normal 2-sphere, we obtain
a cellular decomposition of our type C piece. 
An important observation at this point is that an
``octagonal'' almost normal 2-sphere (with respect
to the original triangulation
$T$)
contained in our
piece, will also be almost normal in this cellular decomposition. 
The details of this construction can again be found
in~\cite{king}.

Like in the case of normal surfaces, the almost normal theory
we outlined above also generalizes naturally to the cellular
setting in our type C piece. Lemma~\ref{lem:t2} gives us
an ``octagonal'' almost normal 
2-sphere in the cellular decomposition of our piece. 
This 2-sphere can be expressed as a sum of the fundamental 
surfaces in the type C piece we are considering. Precisely
one of the summands has to contain a single 
octagonal piece and, since the Euler characteristic is
additive, at least one of the fundamental surfaces in the 
sum has to be a 2-sphere (since the type C piece we are looking 
at contains an ``octagonal'' almost normal 2-sphere, it has to
be a 3-ball and can thus not contain embedded projective planes).
If the fundamental 2-sphere in the sum does not contain an 
octagon, then it has to be normal and therefore parallel to the
boundary of the piece. 
This is a contradiction because we could then isotope it away from
all the other summands by a normal isotopy. Since regular alteration
is defined up to normal isotopy, this would then make the 
sum (i.e. a 2-sphere) disconnected. So we've found an almost
normal 2-sphere in the cellular decomposition of a type
C piece. Such a 2-sphere, looked at from the point of
view of our original triangulation
$T$,
will be an ``octagonal'' almost normal 2-sphere.

The complexity of the fundamental 
``octagonal'' almost normal 2-sphere we've just constructed
is bounded in the same way as all the other complexities of
the normal 2-spheres in
$\Sigma$.
This follows directly from the construction, since all we are
doing when searching for an almost normal 2-sphere, is
just making another step of the recursion that gave us 
$\Sigma$.
We will give an explicit estimate for the complexity 
later on in this section.

%Let
%$\mathcal{A}$
%be a finite system of linear equations such that a 
%non-negative integral solution to 
%$\mathcal{A}$
%corresponds to a (possibly immersed) surface
%$F$, where, except in
%$H$,
%$F$
%consists of normal discs. In 
%$H$,
%$F$
%is composed of normal triangles and octagonal components
%with boundaries parallel to 
%$c$. 
%Just as before, there exists a finite collection 
%$F_1, \ldots ,F_n$
%of fundamental solutions to
%$\mathcal{A}$.
%Again any solution to
%$\mathcal{A}$
%can be written as a finite sum of fundamental solutions.
%At this 
%stage, Thompson~\cite{thompson} proves the following
%lemma.
%
%\begin{lem}
%\label{lem:an}
%If there exists an almost normal 
%2-sphere in the fixed type \textrm{C} piece, then there 
%exists one which is a fundamental solution to some linear
%system 
%$\mathcal{A}$.
%\end{lem}
%
%In order to conduct a full search for almost normal 
%2-spheres in type C pieces, we need to repeat the above
%process for all possible linear systems corresponding
%to all the tetrahedra in 
%$T$
%and
%all the normal curves of length eight on them. 

Let's first bound the number of disjoint non-parallel normal
2-spheres in
$\Sigma$.
This
is made possible by an old idea due to Kneser.  

\begin{lem}
\label{lem:sphere}
Let
$T$
be any triangulation of
$S^3$
and let
$t$
be the number of tetrahedra in 
$T$.
Then any family of disjoint non-parallel normal 2-spheres contains
at most 
$6t$
of them.
\end{lem}

\begin{proof}
Normal triangles and squares chop up any tetrahedron in 
$T$
into several pieces. But at most six of these regions are not of
the form
$\mathrm{triangle}\times I$
or
$\mathrm{square}\times I$
(see figure~\ref{fig:comp}).

\begin{figure}[!hbt]
 \begin{center}
%  \vspace{2cm}
%  \Huge{NO PRODUCT REGIONS}
    \epsfig{file=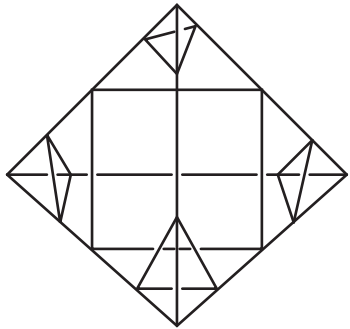}
  \caption{\small Complementary regions which are not a product.}
  \label{fig:comp}
 \end{center}
\end{figure}

Let
$n$
be the maximal number of disjoint non-parallel normal 2-spheres in 
$T$.
Then the complement of this family has precisely
$(n+1)$
components. Each of those components must contain at least one
of the non-product regions. This is because any component, consisting
only of product pieces, is bounded by two parallel normal 2-spheres.
Since the total number of non-product regions is bounded by 
$6t$,
our lemma is proved.
\end{proof}

We are interested in bounding the number of normal pieces of
elements in
$\Sigma$.
We also want to bound
the number of normal 
pieces of the ``octagonal'' almost 
normal 2-spheres that arise
in type C pieces of the cut open 3-manifold
$M-\Sigma$.
Both of these things can be accomplished 
at one go, because we know that the procedure giving
$\Sigma$
can be extended (by making a single additional step) to an
algorithm producing
``octagonal'' almost normal 2-spheres in type C pieces.

The proposition we are about to state is proved
in~\cite{hass}. It originally deals with the linear 
system in
$7t$
variables coming from the matching equations for normal
surfaces. 
Its proof uses some basic linear algebra on the 
linear system which consists of matching 
equations. 
We should note at this point that the number of equations
in this linear system does not influence the bound
that the proposition gives. 

%The same techniques can be applied to the linear 
%system
%$\mathcal{A}$
%we introduced earlier to search for almost normal surfaces.
%The bound from proposition~\ref{prop:hass} will remain 
%valid because 
%$\mathcal{A}$
%is a linear system of the same kind as before (i.e. a homogeneous
%one) and
%it has fewer 
%(precisely
%$7(t-1)+4$)
%variables.

\begin{prop}
\label{prop:hass}
Let 
$M$
be a triangulated 3-manifold containing 
$t$
tetrahedra. 
Let
$\underline{x}$
be a fundamental solution of a system of
linear equations coming from matching conditions. Then
each coordinate of the vector 
$\underline{x}$
is bounded above by
$4 t2^{7t}$.
\end{prop}

Using proposition~\ref{prop:hass}, we can bound the size of
each component of all the vectors corresponding to the 
normal 2-spheres in
$\Sigma$.
It follows from figure~\ref{fig:comp} that the number of 
non-trivial cellular regions (in the cut
open manifold) at any stage of the algorithm
producing the family 
$\Sigma$,
is always bounded by
$6t$.
Proposition~\ref{prop:hass} implies that there can be
at most
$4\cdot 6t\cdot 2^{42t}$
parallel copies of a given normal disc type in a
cellular structure coming from any stage of the algorithm. 

Lemma~\ref{lem:sphere} tells us that we'll never have to
make more then
$6t$
steps when constructing
$\Sigma$.
This means that each of the normal disc types
in the relative cellular structures can only give rise to
less then
$(4\cdot 6t\cdot 2^{42t})^{6t}$
normal discs of the same type in the initial 
triangulation
$T$.

We can obtain a similar kind of bound for 
``octagonal'' almost normal 2-spheres. We only have
to change the exponent from
$6t$
to
$(6t+1)$.
This is because all these ``octagonal'' almost normal 
2-spheres are just one step away (in our algorithm)
from the normal ones 
(bounding type C pieces) and at each stage they are 
described by fewer variables. For example, in our
original triangulation they require 
$7(t-1)+4$
variables. 
So the bound in proposition~\ref{prop:hass} applies.

Putting everything together, we get the following
lemma.

\begin{lem}
\label{lem:bound}
Let 
$T$
be a triangulation of the 3-sphere which contains
$t$
tetrahedra.
Then the number of all normal pieces contained both in all
elements of 
$\Sigma$
and in all ``octagonal'' almost normal 2-spheres 
from all type C pieces is bounded above by
$2^{300t^2}$.
\end{lem}

We should note at this point that this is
the only part of the bound in theorem~\ref{thm:main}
which contains a quadratic expression in its
exponent. If one could find both the
``octagonal'' almost normal 2-spheres in type C pieces
and the maximal family 
$\Sigma$
among the fundamental solutions
of linear systems that are based on the
triangulation
$T$, 
the bound
in lemma~\ref{lem:bound} would have a linear function
(similar to the one in proposition~\ref{prop:hass})
in its exponent.\\ 

The essential process we are just about to describe, 
is the one of 
isotoping almost normal surfaces. It is going to provide
a foundation for the simplifying procedure needed for
the proof of theorem~\ref{thm:main}.

Let 
$F$
be a separating almost normal surface in a 3-manifold
$M$
with a triangulation
$T$.
Its 
\textit{weight},
$w(F)$,
is defined to be the number of points in the intersection 
of
$F$
and the 1-skeleton 
$T^1$.
If 
$F$
contains an octagon,
a
\textit{natural isotopy}
is the one pushing the surface over an edge which meets
the length eight normal curve bounding the octagon in
two points. There are two
possible natural isotopies, depending on the component of
$M-F$
we are pushing into. In case of other non-normal pieces (see
figure~\ref{fig:almostnormal}), a natural isotopy pushes the tube
part of the annulus so that it encompasses one of the edges it is
parallel to. 
As a result 
in both cases, we get a surface with its weight equal to
$w(F)-2$. 

%\begin{figure}[!hbt]
% \begin{center}
%  \vspace{2cm}
%  \Huge{NATURAL ISOTOPY}
%    \epsfig{file=ni.eps}
%  \caption{\small Three possibilities after a natural isotopy over
%                  the edge $f$.}
%  \label{fig:ni}
% \end{center}
%\end{figure}

Notice that if we look at our almost normal surface
$F$
in the complement of the 1-skeleton of
$T$,
there is an obvious compression disc 
$D$
for it, enveloping the edge
we are isotoping over.
The natural isotopy can then be realized by isotoping over the
3-ball bounded by 
$D$
and the disc in
$F$,
bounded by 
$\partial D$.

The natural isotopy is only the first step in the process
of isotoping almost normal surfaces. Everything else will 
be accomplished by a sequence of elementary isotopies. We can define 
them as follows.

\begin{figure}[!hbt]
 \begin{center}
%  \vspace{2cm}
%  \Huge{NON-NORMAL ARC}
    \epsfig{file=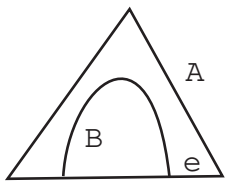}
  \caption{\small An isotope of $F$ intersects $A$ in a non-normal arc.}
  \label{fig:non-normal}
 \end{center}
\end{figure}

Let
$A$
be a 2-simplex in
$T$,
containing a non-normal arc (figure~\ref{fig:non-normal})
which comes from intersecting
$A$
with an isotope of
$F$.
\noindent A disc
$B$
in the triangle
$A$
(see figure~\ref{fig:non-normal})
is bounded by the non-normal arc and a subarc of the edge
$e$.
An
\textit{elementary disc}
can be constructed by banding together two parallel copies
of 
$B$
in the complement of 1-skeleton, where the 
band runs around the edge 
$e$. 
Its boundary is a simple
closed curve in the surface, bounding a disc on one side.
An
\textit{elementary isotopy}
is an isotopy over the 3-ball bounded by the disc in the isotope
of 
$F$
and the elementary disc we've just defined.

Since
$F$
is a separating surface, we can fix a complementary component
$I$
of
$M-F$. 
All the elementary isotopies that we are going to do from now on,
are going to have the same direction. We will always be
isotoping towards the interior of the component
$I$.

The following 
isotoping lemma will play a crucial role in the simplifying
process. A similar result is proved in~\cite{stocking} by
a careful inspection of all the possible cases. 
The proof we are giving here is based on elementary isotopies 
and is better suited from our perspective because it
sheds more light
on the side of things we'll be interested in later.

\begin{figure}[!hbt]
 \begin{center}
%  \vspace{2cm}
%  \Huge{THE END OF ISOTOPY}
    \epsfig{file=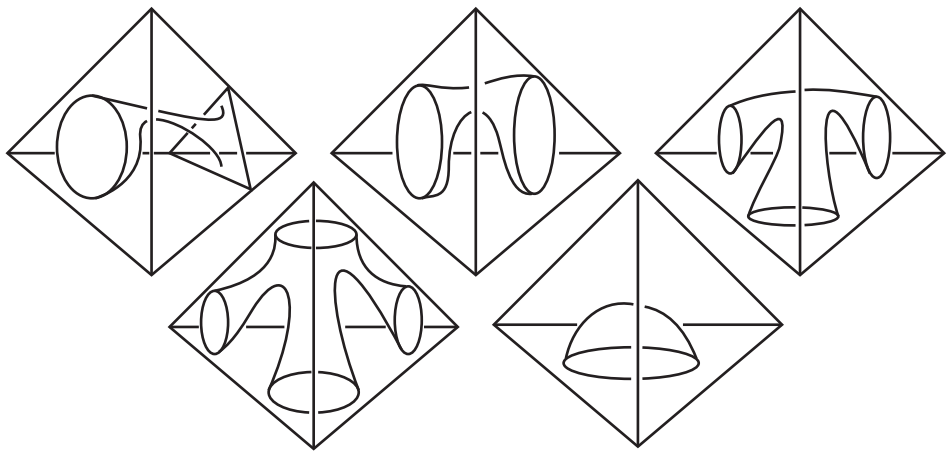}
  \caption{\small Non-normal pieces in the tetrahedra of $T$.}
  \label{fig:end}
 \end{center}
\end{figure}

\begin{lem}
\label{lem:iso}
Let
$F$
be a separating almost normal surface in a 3-manifold
$M$
with a triangulation
$T$.
Let
$I$
be a component
of the complement
$M-F$,
if 
$F$
contains an octagon.
Otherwise let
$I$
denote the component containing a solid torus region 
in the interior of the 3-simplex where 
$F$
is not normal.
A natural isotopy followed by a sequence of 
all possible elementary isotopies,
all going in the direction of 
$I$,
will result with a surface intersecting each tetrahedron
of
$T$
in pieces as in figure~\ref{fig:end} and in 
normal pieces.
Moreover, 
in each tetrahedron
there can only be at most one piece
of the first type from figure~\ref{fig:end}.
A single 3-simplex can contain
several pieces of all the 
other types in figure~\ref{fig:end} as well as several 
normal pieces. The pieces in figure~\ref{fig:end} 
can not be parallel.
\end{lem}

\begin{proof}
First note that after the natural isotopy, all the 
non-normal arcs we get will give rise to elementary isotopies
in the direction of
$I$. 
After each isotopy both 
$F$
and
$I$
will change,
but we'll still denote both resulting spaces by
$F$
and
$I$
respectively.

After the natural isotopy, 
$F$
and
$I$
satisfy the following conditions:
\begin{enumerate}
\item In each tetrahedron of
      $T$
      the component 
      $I$
      consists of a family of 3-balls, each one bounded by pieces of
      $F$
      and a (possibly disconnected)
      planar surface, contained in the boundary of the 3-simplex.

\item Each 3-ball from 1 intersects any face of the tetrahedron it
      lies in, in at most one disc.
\end{enumerate}

An elementary isotopy moves a disc in
$F$
over a 3-ball in
$I$,
which intersects a single edge 
$e$
in
$T^1$.
So the new
$I$
is just the old 
$I$
without the 3-ball we isotoped over. This 
3-ball is a union of a family of 3-balls, one in each 
tetrahedron of the star of
$e$.
In fact there can be more than one 3-ball from the same family
in a single 3-simplex if this 3-simplex occurs
more than once in the star of the edge
$e$.
This is perfectly feasible in a 
non-combinatorial triangulation, but 
it does not have any effect on the process we are studying.

The elements of the above family are the ones that are going
to determine the topology of the pieces of
$I$
in tetrahedra
of
$T$.
In fact, each 3-ball from condition 1 will after an elementary isotopy
still satisfy both conditions if we substitute the old 
$F$
with the new one.
So after performing all possible elementary isotopies towards the
interior of
$I$,
the surface 
$F$
we end up with will intersect each triangle of
$T$ in normal arcs and simple closed curves which miss the
boundary of the triangle.

\begin{figure}[!hbt]
 \begin{center}
%  \vspace{2cm}
%  \Huge{DURING THE ISOTOPY}
    \epsfig{file=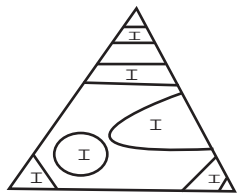}
  \caption{\small An intermediate state of the isotopy on a triangle in $T$.}
  \label{fig:intermediary}
 \end{center}
\end{figure}

The region
$I$
will,
after the isotopy, consist of 3-balls in each 3-simplex.
There is going to be a bijective correspondence between the 3-balls
in the end, and the ones we started with. By condition 2, 
every 3-ball will still intersect any face of the 3-simplex it lies
in, in at most one disc. It is also true that the number of these
discs will not increase when we pass from the 3-ball pieces
of
$I$
at the beginning to the 3-balls at the end.

Let's look at the pieces of
$F$
in each tetrahedron. It is obvious that all the possibilities
of the lemma can actually arise. We have to see that they are the 
only ones. \\
\textbf{Claim.} A single piece of 
$F$
can intersect a triangle of
$T$
in either a unique normal arc or in a single simple closed curve.\\
\indent Every piece of 
$F$
is contained in the boundary of a 3-ball piece of
$I$. 
This 3-ball intersects each triangle of the 2-skeleton 
$T^2$
in at most one disc. So no triangle can contain two simple closed
curves or a simple closed curve and a normal arc, both belonging 
to the same piece of
$F$.

The same argument tells us that a triangle
in 
$T^2$
can either contain 
two normal arcs
of intersection with a single piece of 
$F$
or at most three of them, each one cutting
off a vertex of the triangle
in the 2-skeleton. 

Now we need to prove that our piece of
$F$
can have at most one normal 
simple closed curve boundary component.
So assume the opposite.
Since the piece is a subset of the boundary of
a 3-ball, no arc contained in it, running between 
two distinct boundary components
of our piece, can be extended to a simple 
closed curve in the 2-sphere bounding that 3-ball, without 
increasing 
the number of intersection points with the boundary of our piece.
On the other hand,
assuming we have at least
two normal simple closed curves in
the boundary, there surely exist two normal arcs,
belonging to the distinct boundary components of our piece,
that are contained in a single 2-simplex.  
Connecting them by an arc in the piece of
$F$
contradicts what was said before (because these two normal
arcs are both contained in the boundary of a disc in the 2-simplex
they lie in).

So now it follows 
that the piece of
$F$
we are looking at, can contain at most 
one normal boundary component which is of length at most eight.
This is because the only normal curve
of length 12, intersecting each 2-simplex (in the boundary
of a tetrahedron) in 3 normal arcs as above, consists of
4 simple closed curves (one 
for each vertex of the tetrahedron). It is also well known that 
normal simple closed curves of lengths 9, 10 or 11 do 
not exist.

There are precisely three normal simple closed 
curves of length eight in the boundary of
a tetrahedron.
So if our piece of 
$F$
is bounded by one
such curve,
then
at least one of the faces of the tetrahedron intersects
the 3-ball piece of
$I$
(containing in its boundary the piece of 
$F$
we are considering)
in two discs. This is a contradiction that proves the claim.

The claim implies the following seven possible boundaries 
for any piece
of 
$F$: normal simple closed curve of length three, normal 
simple closed curve of length four, single simple closed curve,
normal simple closed curve of length three and a 
simple closed curve, two 
simple closed curves,
three simple closed curves, four simple closed curves.

Since every piece of 
$F$
in any tetrahedron is planar, it is 
up to homeomorphism determined by its boundary.
This implies that all possible pieces of 
$F$
are the ones listed in the lemma.

The fact that all these planar surfaces are embedded 
as in figure~\ref{fig:end} 
(up to an isomorphism of the tetrahedron) follows from the
observation that all the elementary discs are parallel 
to edges of the 1-skeleton. 
\end{proof}

\begin{center}
\section{\normalsize \scshape OUTLINE OF THE PROOF}
\label{sec:out}
\end{center}

Given a triangulation
$T$
of the 3-sphere,
how do we simplify it? The process is divided into two stages. First, we 
create a subdivision 
$S$
of 
$T$
by defining it in each complementary piece of the manifold
$S^3-\Sigma$
in such a way that the triangulations match along
all normal 2-spheres in
$\Sigma$.
The second step consists of simplifying 
$S$
down to the canonical triangulation of 
$S^3$.

An explicit construction of
$S$,
using Pachner moves,
will be given in section~\ref{sec:sub}. The simplifying
procedure 
of step two is based on the relationship between Pachner
moves and shellable triangulations. 
This relationship will be established
in section~\ref{sec:pach}. 

Now, we are going to describe the additional structure 
on the complementary pieces of
$S^3-\Sigma$,
needed for the definition of the subdivision
$S$.
We already know (lemma~\ref{lem:t2}) that every type C
piece contains an ``octagonal'' almost normal 2-sphere. 
To see that each type B piece also contains an almost
normal 2-sphere, it is useful
to introduce an ordering on the normal family
$\Sigma$.
It comes
naturally by picking a vertex of
$T$
and looking at the complementary region 
(which is not a type A piece) of the trivial normal 2-sphere 
around it.
Topologically we get a 3-ball containing our normal family
$\Sigma$.
Now, the ordering on
$\Sigma$
is induced by inclusion. For example, the trivial normal 2-sphere around the
vertex we removed is the unique largest element.
The smallest elements in this ordering are either trivial normal
2-spheres consisting of normal triangles only 
or the ones bounding type C 3-balls.

Our task is to find an almost normal 2-sphere in a
piece with more than one boundary component.
Pick
the largest 2-sphere in its boundary. A very nice argument
in~\cite{thompson} (subclaim 2.0.1.) implies that there must
be an edge in
$T$
with a subarc
which runs from the largest component of the boundary to some other
component
and whose interior is disjoint from 
$\Sigma$.
By taking parallel copies of the two 2-spheres connected by this
arc in the piece we are looking at, and tubing them together
in one of the tetrahedra in the star of the edge, we obtain
our almost normal 2-sphere.

\begin{figure}[!hbt]
 \begin{center}
%  \vspace{2cm}
%  \Huge{TYPE C PIECE}
   \epsfig{file=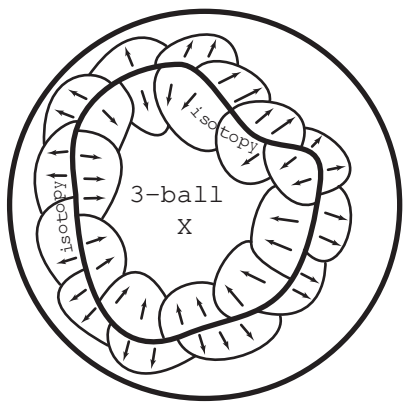}
  \caption{\small Two sequences of elementary isotopies
                  in a type C piece.}
  \label{fig:typeC}
 \end{center}
\end{figure}

\begin{figure}[!hbt]
 \begin{center}
%  \vspace{2cm}
%  \Huge{TYPE C PIECE}
   \epsfig{file=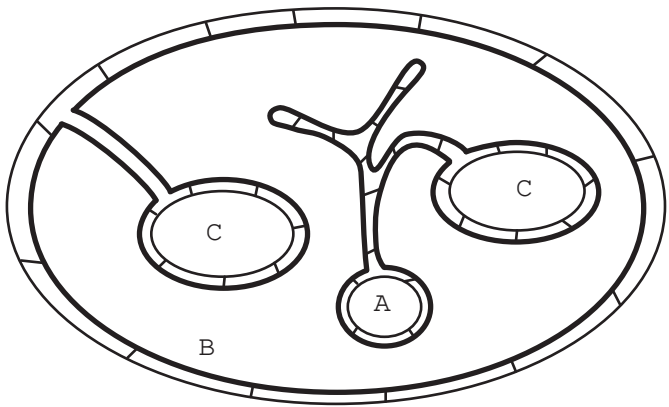}
  \caption{\small The ``tubed'' almost normal 2-sphere and the isotoped 
                  2-sphere in a type B piece.}
  \label{fig:typeB}
 \end{center}
\end{figure}

All the almost normal 2-spheres we've created are separating,
because we are in
$S^3$.
By picking the right complementary component 
in
$S^3$
and applying 
lemma~\ref{lem:iso}, we can simplify each almost normal
2-sphere by a sequence of elementary isotopies. Since we 
are only using elementary isotopies going in the same
direction (towards the interior of a fixed complementary 
component
in
$S^3$), the whole process can be realized by an
embedding of
$S^2\times I$,
where the top 2-sphere is the almost normal 2-sphere we started
with and the bottom one is the 2-sphere coming from
lemma~\ref{lem:iso}
(see figures~\ref{fig:typeC} and~\ref{fig:typeB}).

Another important point here is that the whole isotopy
never leaves the type B (or C) piece it started in.
This is true simply because an analogous statement holds
for each elementary isotopy.
This implies that the isotoped surface coming from
lemma~\ref{lem:iso} will be contained in the interior
of the piece containing the almost normal 2-sphere we started 
with. 

In a type C piece, the isotopy can go in two directions
because the almost normal 2-sphere in this case contains
an octagon. The surface we get, when isotoping 
towards the interior of the piece, will have 0 weight. 
This follows from the observation that we can forget
about all pieces in figure~\ref{fig:end}
if we compress each annulus with a length three normal curve
in its boundary. This would then give a family of normal 2-spheres
contained in the type C piece 
which is a contradiction. Therefore, the 2-sphere we end up with
has to miss the 1-skeleton.

Similar reasoning tells us that an isotopy in the other 
direction in the type C piece has to end with a 2-sphere,
intersecting the 2-skeleton
$T^2$
in normal curves parallel to the ones coming from 
the boundary of the piece we are looking at and possibly in 
some simple closed curves which miss the 1-skeleton.

The almost normal 2-spheres in type B pieces 
that we are going to consider, will never contain
an octagon. We will therefore be isotoping in one direction
only. Using the same kind of arguments as before, we can
conclude that the 2-sphere we end up with consists of boundary
components
of the type B piece 
(all except the two we started with, 
see figure~\ref{fig:typeB}), tubed together by
pieces depicted in figure~\ref{fig:end}. It should be noted that
if a type B piece has only two boundary components, then the
isotoped 2-sphere does not intersect the 1-skeleton.

Let
$\Lambda$
be
the following collection of 2-spheres: in every type C piece
just take an 
``octagonal'' almost normal 2-spheres which
exists by lemma~\ref{lem:t2}. In each type B piece take
a copy of the almost normal 2-sphere described above
with the annulus connecting two normal pieces moved by a 
natural isotopy,
so that it envelops the edge it is
parallel to. The 2-spheres from 
$\Lambda$
in type B pieces are therefore normal in all the tetrahedra of
$T$,
except in the ones contained in the star of the edge
we isotoped over. 

The sequences of elementary 3-balls, 
corresponding to the supports of elementary isotopies, 
yields the additional structure (on type B and C pieces)
that is required to define
the subdivision 
$S$.
Elementary discs  
and an element in
$\Lambda$
chop up each type C piece of
$S^3-\Sigma$ (see figure~\ref{fig:typeC}).

In the case of a type B piece,
the element of
$\Lambda$
will, after the isotopy, consist of all
but two
boundary normal 2-spheres tubed together by pieces described in the
figure~\ref{fig:end}. 
Again, the type B piece in question can be decomposed into 3-balls
by all the elementary discs required for the
isotopy
and by the element 
in
$\Lambda$
we started our isotopy on. 

After the isotopy from lemma~\ref{lem:iso}, what's left in each
tetrahedron of the complementary component we isotoped into, are
just 3-balls bounded by the pieces from 
lemma~\ref{lem:iso} on one side and possibly some normal pieces of the
elements in
$\Sigma$
on the other. Schematically, the situation after the isotopy 
is depicted by figure~\ref{fig:typeB}.

Now we want to triangulate all 
of these 3-balls (the elementary ones as well as the
ones that are left over in the component we were
isotoping into), in all the pieces of the complement of
$\Sigma$,
by 
simple shellable triangulations.
Since all the processes described above induce
polyhedral structures in the boundaries of all the 
3-balls (this will be described
in detail in section~\ref{sec:sub}) in question,
subdividing the boundary 
2-spheres 
in order to obtain  genuine 
triangulations and then coning them, does the job.
Doing so in every piece of the space
$S^3-\Sigma$
exhausts the whole 3-sphere and therefore completely 
determines the subdivision 
$S$.

The fact that all these cones are indeed shellable, is 
proved in~\cite{lickorish} (lemma 5.4). 
Here we are relying on the property that all the bounding 2-sphere
we'll need to cone in the process, are triangulated by 
combinatorial triangulations. 
The reason why we want
these 3-balls to be shellable is simply because 
each elementary isotopy can then be
realized by a shelling of the corresponding 3-ball.

So what we really want from the cones on the 2-spheres above,
is to be shellable without ever having to shell from the 
faces contained in a fixed disc, which is lying 
in the bounding 2-sphere. 
This disc is just a 2-manifold along which the 3-ball, we are
trying to triangulate, is glued onto the rest of the (type B 
or C) piece. This can always be achieved since a cone on a
disc with a combinatorial triangulation can be shelled
``from the side'' just by coning the shelling procedure
of the disc itself.  

The simplifying process works its way up the ordering of the
normal 2-spheres in 
$\Sigma$.
First we change the subdivision 
$S$
in all type C pieces (which are smallest elements in our ordering), 
making it a cone on the unique boundary component 
in each of the pieces. 
In section~\ref{sec:pach} we will discuss
how to implement elementary shellings from a 2-sphere
boundary component 
by Pachner moves, if on the other side of that 2-sphere
we have a cone on it. Using that construction we can
pick a 3-ball piece 
(in some 3-simplex), 
contained in the
3-ball
$X$ 
from figure~\ref{fig:typeC}, 
and turn the whole 3-ball 
$X$, bounded by the 2-sphere coming out 
of lemma~\ref{lem:iso}, into a cone on its boundary.
This is simply because the complement (in 
$X$) 
of the coned 3-ball
piece we picked, is shellable. 
That follows from the observation that all the 
(coned) 3-ball regions from figure~\ref{fig:end},
the 3-ball
$X$
is made of (we already know
that the first possibility in figure~\ref{fig:end} can 
not occur) 
can be viewed as vertices of a graph whose edges 
correspond to the discs in the interiors 
of the 2-simplices of
$T$.
Since this graph is a tree
(this follows from the fact
that the isotoped 2-sphere bounds 
a 3-ball), there is a
``global'' shelling strategy for the 
complement of the piece we picked in the 3-ball
$X$.
This can be made simplicial by shelling one cone at a time.

So now we can assume that the 3-ball 
$X$
is coned.
We can carry on by shelling (in the reversed order)
all elementary 3-balls (and the
3-ball corresponding to the natural
isotopy) involved in the isotopy taking
the almost normal 2-sphere to the boundary of
$X$. 
By this stage, we've changed the subdivision 
$S$
so that it looks like a cone on the almost normal 
2-sphere in the type C piece we are looking at.
Above it,
$S$
is still unchanged.
We can now do the same thing towards the boundary of the 
piece we are considering, again using the shellable 
nature of the subdivision 
$S$
in all appropriate 3-balls.

What we have now is a cone on the 2-sphere which intersects
the 2-skeleton
$T^2$
in normal curves, parallel to the ones coming from the 
bounding normal 2-sphere of the type C piece we are looking at. 
First we shell all the 3-balls from the 
$S^2\times I$
region, bounded by the isotoped 
``octagonal'' almost normal 2-sphere and the single boundary component 
of the type C piece, that are bounded by pieces in 
figure~\ref{fig:end}. We thus obtain a cone on a normal 2-sphere 
which is parallel to the bounding 2-sphere of our piece. 
Since all the regions between 
any two parallel normal pieces are cones as well, we can shell 
them one by one and therefore get a cone on the boundary of
our piece. Here we are relaying on the fact that the 
normal structure on the bounding 2-sphere is shellable. 
In general this needn't hold, but the technical
assumption that we are going to make on our triangulation
$T$
at the beginning of section~\ref{sec:sub} will guarantee
this property.
This completes the simplification
of the triangulation 
$S$
in all type C pieces.

Take a type B piece and assume that all normal 2-spheres, 
strictly smaller than the largest normal 2-sphere in its boundary,
already bound coned 3-balls. The strategy now is similar to the one
we used in type C pieces. 
Since the tube of the 2-sphere element of
$\Lambda$
in our piece runs
from the largest  boundary component to some other
boundary component, we can deduce that all other 
normal 2-spheres in the boundary that are going to be tubed together
by the isotoped 2-sphere
(see figure~\ref{fig:typeB}), are going to bound 
cones on one side. 

We first shell all the regions
which are bounded by two parallel normal pieces and lie
between a normal boundary component and the isotoped
2-sphere. 
We can do that by expanding the cone structures 
on the other side of the boundary components of
our type B piece, that 
exist by assumption. Now the 3-ball bounded by the isotoped
2-sphere from 
$\Lambda$ 
is again chopped up into 3-ball regions
that are glued together along discs contained in the interiors
of 2-simplices of
$T$.
Like before, there is a sequence of elementary shellings
which gives a way of changing the triangulation of the
3-ball bounded by the isotoped 2-sphere 
from
$\Lambda$
to the 
cone on its boundary.  

We will now mimic what we did in type C pieces. Let's take
the sequence (in the reversed order) of all elementary 
3-balls coming from the elementary isotopies needed to push
the 2-sphere element of 
$\Lambda$
in our type B piece, down to the 2-sphere which now already 
bounds a coned 3-ball. Using this sequence in the same way as above,
we can change the triangulation
$S$
in our piece to the cone on the 2-sphere element of 
$\Lambda$.
It is now obvious how to simplify the remains of the 
subdivision 
$S$
in the type B piece we are looking at. 
So we've managed to transform the subdivision 
in our piece into a cone on the 
largest boundary component. 

Now we want to make sure that the techniques 
described above suffice for the total simplification
of the subdivision 
$S$.
Let's assume that 
$K$
is a triangulation of 
$S^3$
and that
$S$
is its subdivision containing 
$\Sigma$
as a subcomplex. Let 
$x$
be the vertex of
$K$
inducing an ordering on
$\Sigma$.
Let's also
assume the following property: 
if all elements of
$\Sigma$,
smaller then a given normal 2-sphere
$A$
in
$\Sigma$,
bound coned 3-balls, 
then using Pachner moves we can change the
triangulation of the 3-ball bounded by 
$A$
into a cone on 
$A$,
without altering the simplicial structure of
$A$.
Then we claim that we can transform 
$S$,
using Pachner moves only, into a cone
on
$x$
glued to another copy of itself via an identity on the boundary.

To see this, we'll use a simple induction 
on the depth of elements in
$\Sigma$.
A normal 2-sphere in
$\Sigma$
is of \textit{depth} 
$k$
if it is greater than precisely 
$k$
elements of
$\Sigma$.

We can use our assumption for the 2-spheres of
depth 0. 
Let 
$A$
in 
$\Sigma$
be of depth 
$(k+1)$
and assume we've coned all the 2-spheres
of depth smaller or equal to 
$k$.
Any 2-sphere smaller than 
$A$
is of depth at most 
$k$. 
So we can use the assumption again.
This proves our claim and therefore
completes the simplification process.

\begin{center}
\end{center}

\begin{center}
\section{\normalsize \scshape PACHNER MOVES AND SHELLABLE TRIANGULATIONS}
\label{sec:pach}
\end{center}

In this section we are going to establish a relationship
between elementary shellings and Pachner moves. We will do this
in dimensions two and three. Both cases will play a crucial role
in building and simplifying the subdivided triangulation
$S$.
Let's start by stating precisely what we mean by shelling.\\

\begin{defin}
Suppose that 
$M'$
is a submanifold
of a triangulated 
$n$-manifold
$M$
with boundary.
If there exists an
$n$-simplex
$\Delta$
in the triangulation of 
$M$
with the property that 
$\Delta\cap\partial M$
is a combinatorial 
$(n-1)$-disc, 
such that 
$M'$
equals the closure 
(in
$M$)
of the complement
$M-\Delta$,
then we say that 
$M'$
is obtained from 
$M$
by an 
\textit{elementary shelling}.
\end{defin}\\

An elementary shelling is quite similar to an elementary collapse
of the top dimensional simplex. The crucial difference lies in
the fact that here we stipulate explicitly that the resulting
space has to be a manifold. 

Another thing which is worth mentioning is that the boundaries
$\partial M$
and
$\partial M'$
differ by a single 
$n-1$ 
dimensional Pachner move. 

A sequence of such elementary shellings is called a
\textit{shelling}. Saying that a triangulation of
an
$n$-manifold is \textit{shellable} simply means that
there exists a sequence of elementary shellings which will
reduce the triangulation down to a single 
$n$-simplex.
Since the homeomorphism type of the manifold in question 
does not change under an elementary shelling, it is clear
that 
$n$-balls 
are the only candidates to have shellable triangulations.
It is for example very well known that any 
combinatorial triangulation
of the two dimensional disc is always shellable. As it was
mentioned before, the lemma 5.4 in~\cite{lickorish} and the
above observation about discs together imply that a cone on
any combinatorial triangulation of 
the 2-sphere constitutes a shellable triangulation of the 3-ball.

Now we are going to express all possible
elementary shellings by Pachner moves in the following 
three dimensional
situation. 
Suppose we had a triangulated 3-manifold and we wanted to make
an elementary shelling from a 2-sphere boundary component. Suppose
further that on the other side of this 2-sphere, we had a cone on it.
We have to consider three different cases according to the
number of faces of the 3-simplex we are shelling, which are contained in the
boundary 2-sphere. 

\begin{figure}[!hbt]
   \begin{center}
%     \vspace{2cm}
%     \Huge{1 FACE}
    \epsfig{file=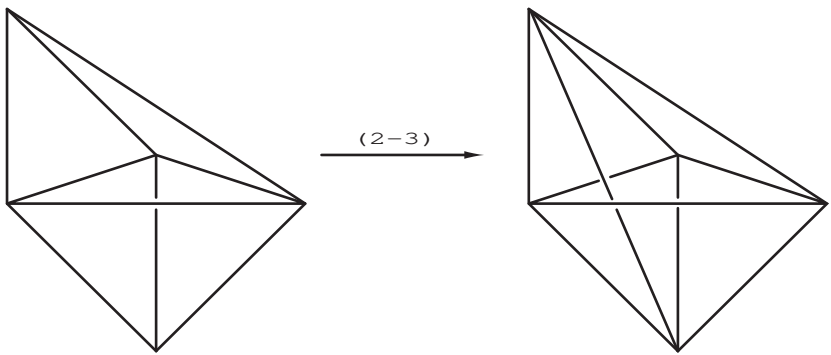}
        \caption{{\small A single free face requires one $\pt$ Pachner move.}}
        \label{fig:1f1}
   \end{center}
\end{figure}

The first case, where we have a single triangle in the boundary,
is dealt with by figure~\ref{fig:1f1}.

We should note that before making the 
$\pt$
move in figure~\ref{fig:1f1}, the top 3-simplex is
contained in the manifold, while the bottom one
belongs to the cone. After the move, all three 3-simplices
are contained in the altered cone.

\begin{figure}[!hbt]
   \begin{center}
%     \vspace{2cm}
%     \Huge{2 FACE}
    \epsfig{file=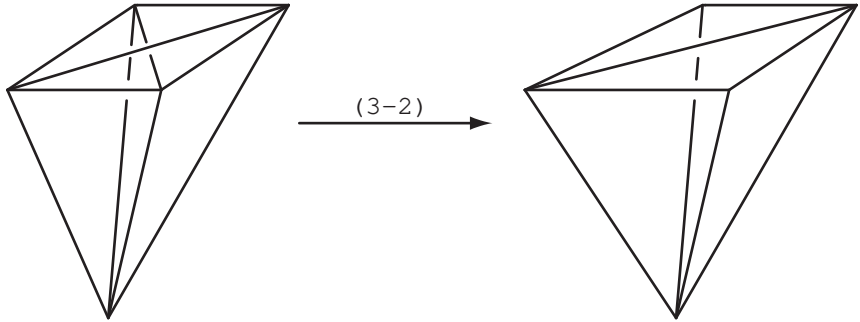}
        \caption{{\small A single $\pth$ move implements the shelling
                         with two free faces in the boundary.}}
        \label{fig:3-2}
   \end{center}
\end{figure}

The second case is the one where we have two faces in the boundary.
It is clear from figure~\ref{fig:3-2}, that a single
$\pth$
Pachner move suffices.

\begin{figure}[!hbt]
   \begin{center}
%     \vspace{2cm}
%     \Huge{3 FACE}
    \epsfig{file=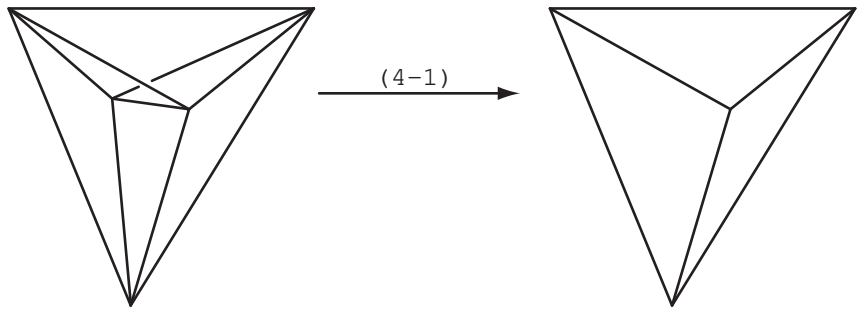}
        \caption{{\small A single $\pf$ move completes the 
                         elementary shelling.}}
        \label{fig:4-1}
   \end{center}
\end{figure}

Finally, we have to deal with the situation where the 3-simplex
we want to shell has three of its faces in the boundary 2-sphere.
The top 3-simplex on the left of figure~\ref{fig:4-1} is the
one we want to shell next, while the other three are contained in the 
cone. It is obvious that a single
$\pf$
Pachner move does the job.

Putting all these facts together, we've seen that in the 
setup described above, each elementary shelling corresponds 
to a single Pachner move. So if we want to bound
the number of Pachner moves required for the simplification
of the subdivision
$S$,
all we need to do is to count the number of tetrahedra in
$S$.
This will be dealt with in section~\ref{sec:conc}.

Before we go on to discuss the two dimensional case, we 
need to prove the following slightly technical lemma
which connects collapsing of an edge with Pachner moves.
It will be of use to us in section~\ref{sec:sub}.

\begin{lem}
\label{lem:crush}
Let
$x$
be a vertex in a combinatorial triangulation of
$S^2$
containing
$n$
2-simplices. Assume further that the star of
$x$
is an embedded PL disc, triangulated by
$k$
triangles. Let
$e$
be the unique edge in the
3-ball, triangulated as a cone on
$S^2$,
running between
$x$
and the cone point.
The triangulation of the same 3-ball obtained by
crushing the edge
$e$,
and thus flattening its star,
can be constructed by 
$(n-k+1)$
Pachner moves used on the original (coned) triangulation.
\end{lem}

\begin{proof}
The 3-ball from the lemma can be view as a union of the following
two PL 3-balls: the star of the edge
$e$
and the cone on the disc in the bounding
$S^2$, which is the complement of the star of
the vertex
$x$
on the 2-sphere.

The triangulation we are aiming for is equal to the triangulation
of the latter 3-ball. We therefore want to flatten the star of the
edge
$e$
down to the cone on the link of
$e$.

This can be achieved by ``moving'' the cone points of the
3-simplices in the second of the two 3-balls described above,
from our initial cone point to the vertex
$x$.
Such a 3-simplex, having a face in
$S^2$,
which is adjacent to the star of
$x$,
can be moved by a
$\pt$
move
or its inverse,
depending on the number of edges it has in common with the
star of
$x$.

Repeating this for all (but one) 3-simplices in the cone on the
disc
$S^2-\i(\mathrm{star}(x))$
almost does the job. All we have to do at this stage, is to use a single
$\pf$
move on what's left of the two 3-balls described above.

We should also note that the sequence of
$\pt$
moves and their inverses, we used to alter
the initial triangulation, can always be found. This follows from
the well-known fact that every combinatorial triangulation of a
PL 2-disc is shellable.
\end{proof}

The rest of this section will be devoted to two dimensional
Pachner moves and their relationship with elementary shellings.
In fact, what we want to do is to transform any given 
triangulation of a disc 
into a cone on its boundary, using 
Pachner moves only. 

\begin{figure}[!hbt]
  \begin{center}
%     \vspace{2cm}
%     \Huge{(2-2) (1-3)}
    \epsfig{file=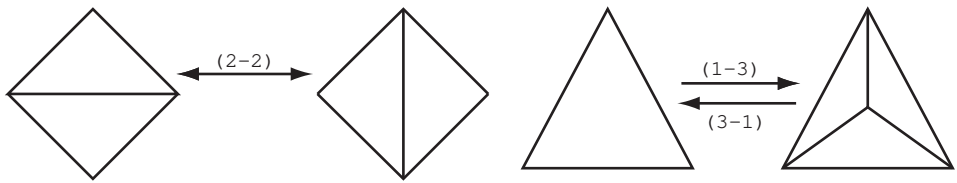}
        \caption{{\small Two dimensional Pachner moves.}}
        \label{fig:2pm}
   \end{center}
\end{figure}

In dimension two
there are three possible moves at our disposal. They are given by
figure~\ref{fig:2pm}.

The simplifying procedure for any PL disc is described by the next
lemma.

\begin{lem}
\label{lem:disc}
Any combinatorial triangulation of a piecewise linear disc
with 
$n$
triangles can be altered into a cone on the boundary of the same
disc by 
$n$
Pachner moves.
\end{lem}

There are two reasons why we can 
(and have to)  assume that the triangulation of 
the disc is combinatorial. The first one is that in what follows,
we can easily guarantee this property for all the discs we are going
to be using our lemma~\ref{lem:disc} on. The second one is that 
the proof of the above lemma relies on the fact that any triangulation
of a disc is shellable, a fact not entirely correct (with our 
definition of an elementary shelling) if we allow for non-combinatorial
triangulations.

\begin{proof}
Since the triangulation of our disc is shellable, we can index all
the simplices in it
by numbers from
$1$
to
$n$, 
so that the increasing integers specify a way
of reducing our triangulation down to a single triangle. The 
2-simplex that's left has index
$n$.
Let's make a 
$(1\p 3)$
move on it. The 2-simplex corresponding to
$n-1$
has to share a unique edge with it. Making a 
$(2\p 2)$
move over this edge, changes our original triangulation
in the last two 2-simplices to a cone on the boundary
of the disc that they compose. The rest of the triangulation
is unchanged at this stage.

Noticing that the union of the last 
$k$
2-simplices in our sequence always gives a disc, makes the 
following induction possible. Say that we already have 
a cone on the boundary of the disc which is the
union of the last 
$k$
2-simplices and that the rest of the triangulation we started with
is unchanged. If the triangle corresponding to 
$n-k-1$ 
has a single edge in common with our cone, we act as before 
(a single 
$(2\p 2)$
move suffices).
If it has two faces in common, a single
$(3 \p 1)$
move finishes the proof.
\end{proof}

\begin{center}
\end{center}

\begin{center}
\section{\normalsize \scshape  THE SUBDIVIDED TRIANGULATION}
\label{sec:sub}
\end{center}

Let 
$T$
be a possibly non-combinatorial triangulation of
$S^3$
with
$t$
tetrahedra.
Let's also make the following technical assumption on
$T$:
each edge in the 1-skeleton of 
$T$
appears at most once as an edge of any 3-simplex in
$T$.
This assumption does not imply that the triangulation
$T$
is combinatorial, but it is certainly satisfied by 
all combinatorial
triangulations of
$S^3$. 
We are making it at this stage because it is 
going to simplify some of the processes we'll have to
invoke later on. It will also become clear that any 
triangulation can be altered so that it has this
property by linearly (in
$t$)
many Pachner moves.

In this section, we shall describe the subdivision 
$S$
of the triangulation 
$T$
and also bound the number of Pachner moves required to construct 
it.

Let 
$\Gamma$
be
the union of all discs needed to perform all 
elementary and natural isotopies in all the pieces of
$S^3-\Sigma$. 
We should note that the number of elements of
$\Gamma$,
coming from a single 2-sphere
in
$\Lambda$,
is bounded above by the number of times 
the almost normal sphere in question intersects the
1-skeleton.
An explicit bound on the number of elements of
$\Gamma$
will be given later.

\begin{figure}[!hbt]
  \begin{center}
%     \vspace{2cm}
%     \Huge{REGIONS OF A GLUED ON DISC}
    \epsfig{file=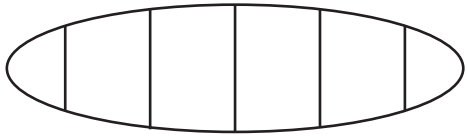}
        \caption{{\small Regions, in a disc in $\Gamma$, correspond to
                         tetrahedra in the star of an edge.}} 
        \label{fig:disc}
   \end{center}
\end{figure}

An elementary disc from
$\Gamma$
will intersect every tetrahedron in the star of the edge we are isotoping
over in a disc region (see figure~\ref{fig:disc}). 

In each tetrahedron, the operation of adding in this disc will consist
of gluing in a length four (respectively two) disc so that two 
(respectively one)
arcs in its boundary are contained in the isotoped almost 
normal surface, and
the other two (respectively one) lie in the boundary of the tetrahedron. 

We are now in the position to describe the subdivision
of the polyhedron
$$T^2\cup\Sigma\cup\Lambda\cup\Gamma$$
which will be a subcomplex of the triangulation
$S$. 
In fact, the simplicial structure of the polyhedron
$T^2\cup\Sigma\cup\Lambda\cup\Gamma$
will play a crucial role in the simplifying process and will also 
be of significance in the definition of the subdivision 
$S$.

All the normal 2-spheres in
$\Sigma$
will inherit the PL structure from their normal structure. 
The normal triangles in 
$\Sigma$
will become 2-simplices, while the normal
quadrilaterals will be subdivided into two 2-simplices
by a diagonal. 

The PL structure of the almost normal 2-spheres in
$\Lambda$
will be a subdivision of the normal and almost normal pieces.
We will subdivide them according to the markings on them, 
made by discs
in
$\Gamma$, 
where we define a
\textit{marking}
on a normal or an almost normal piece to be an arc 
of intersection of the piece with a disc in
$\Gamma$.
We know, that each element of 
$\Gamma$ 
is chopped up into discs of lengths two or four in each tetrahedron
(as in figure~\ref{fig:disc}). 
%One or two of the boundary arcs
%of a region in an element in
%$\Gamma$
%will be embedded in  normal or in  almost normal pieces. 
%This process, in the case of two parallel triangles belonging to the same
%almost normal 2-sphere, is depicted on figure~\ref{fig:marking}.
%
%\begin{figure}[!hbt]
%  \begin{center}
%     \vspace{2cm}
%     \Huge{MARKING}
%    \epsfig{file=protpri2.eps}
%        \caption{{\small Markings of two parallel normal triangles.}}
%        \label{fig:marking}
%   \end{center}
%\end{figure}

The disc regions of the elements of
$\Gamma$
of length four can only leave a marking 
on a normal piece of an almost normal 2-sphere going
from one normal arc to another.
There can only be three such markings on a triangle and four of them
on a quadrilateral,
one for each corner. On an almost normal octagon, immediately after we
glue on our first disc
corresponding to 
the natural isotopy, we end up with two triangles. So there can be
at most six markings on an octagon, coming from the discs of length
four.

\begin{figure}[!hbt]
  \begin{center}
%     \vspace{2cm}
%     \Huge{TRIANGLE QUADRILATERAL OCTAGON}
    \epsfig{file=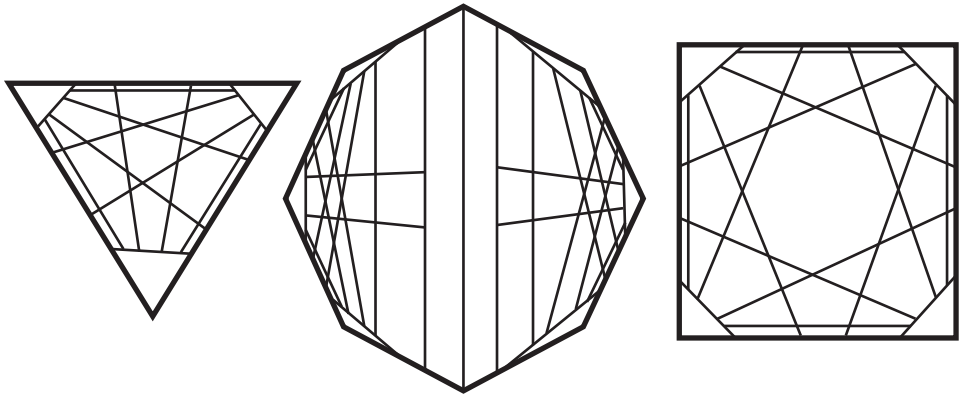}
        \caption{{\small The polyhedral structure of normal and
                         almost normal pieces of elements in $\Lambda$.}}
        \label{fig:superim}
   \end{center}
\end{figure}

The discs of length two will either leave a marking running from
a normal arc to some other 
marking or simply running between two markings.
Because 
each marking is parallel to some edge of the normal piece that it lies
on
and because we can not get more than one marking of the same kind,
superimposing all the possible markings on normal and almost normal
pieces is described by figure~\ref{fig:superim}.

The almost normal piece which is obtained by tubing together two
normal pieces can be treated in the same way, since we could view
it as an annulus around an edge between two normal 2-spheres. 
This annulus consists of discs of length four in each tetrahedron in the
star of the connecting edge. These discs will be glued on the pairs
of normal pieces yielding non-normal pieces, similar to the ones
we get during the isotopy of the surface
$F$
from 
lemma~\ref{lem:iso}.
The PL structure on such a piece will
come from the PL structure on the two parts of normal pieces it
consists of, and from
the PL structure on the glued in disc, which we haven't yet described.
These glued in discs from the annulus behave in the 
same way as the disc regions from elements in 
$\Gamma$
(see figure~\ref{fig:disc}). 

In the next paragraph
we shall see that each of 
these disc regions can be triangulated by at most 6 triangles.
Counting the regions in the normal and almost normal 
pieces in figure~\ref{fig:superim} and triangulating 
each region (if it is not a triangle already)  by coning 
from one of the vertices in its boundary, we can see that 
each piece, including the ones coming from the 
``tubed'' almost normal 2-spheres, contains less than 200
2-simplices. We should also note that the described 
subdivision of the pieces is combinatorial.

Now, we have to put a PL structure on the elements of 
$\Gamma$. We've noted before (figure~\ref{fig:disc}) that each 
elementary disc
in
$\Gamma$
consists of disc regions of lengths two or four. Once we've glued
in a disc from 
$\Gamma$,
the disc regions in it give us a polyhedral structure on it. Further 
gluings will however subdivide this structure. Concentrating on a single
disc region 
$A$
of our element in
$\Gamma$,
we note that all further gluings of disc regions 
of length four will miss 
$A$
completely and therefore not change it at all.
Disc regions of length two can add in a further arc on 
$A$
which runs parallel to the arc(s) in its boundary, contained
in the 2-skeleton
$T^2$. 
Since this can only happen once per boundary arc of 
$A$
in the 2-skeleton
$T^2$,
we can add at most two arcs in each disc region of any 
element in
$\Gamma$.
So a disc in 
$\Gamma$
will in the end look exactly like the disc in figure~\ref{fig:disc}
with less than 
$3t$
disc regions.
This follows from the assumption we made at
the beginning of this section, since it
implies that a star of an edge can contain at most
$t$
tetrahedra.

The arcs in the boundaries of disc regions of elements
in
$\Gamma$
that leave markings on normal and almost normal
pieces of the elements in
$\Lambda$
will be subdivided further by the vertices
coming from the points of intersection of the markings
(see figure~\ref{fig:superim}). An arc in the boundary
of the length two disc region (i.e. the one that's
leftmost or rightmost in figure~\ref{fig:disc}) will
get at most 16 vertices in this way, while an arc in the
boundary of the length four disc region will contain at most
5 such vertices (see figure~\ref{fig:superim}). 

All these observations about the polyhedral structure
of the discs in
$\Gamma$
imply that each disc region corresponding to a single tetrahedron
in
$T$,
will be triangulated by 
no more than 20 triangles. 
So we can triangulate any element from 
$\Gamma$
by less than
$20t$
triangles. Again, the triangulation we get is combinatorial.

\begin{figure}[!hbt]
  \begin{center}
%     \vspace{2cm}
%     \Huge{REGIONS}
    \epsfig{file=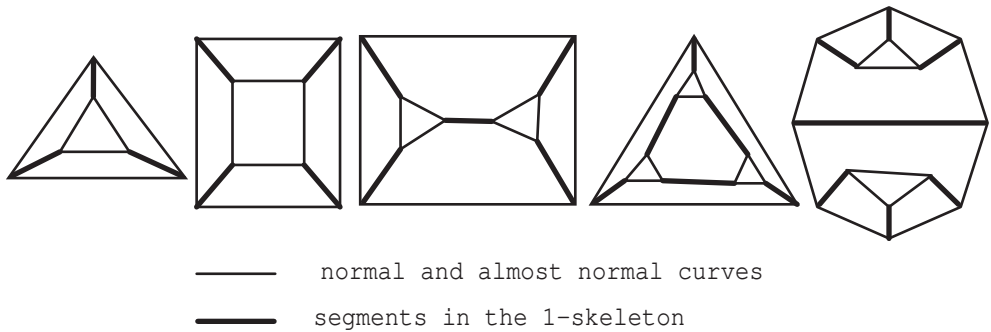}
        \caption{{\small Regions in the boundary of a 3-simplex
                         bounded by normal and almost normal simple
                         closed curves.}}
        \label{fig:regions}
   \end{center}
\end{figure}

Finally, we need to induce a PL structure on the 2-skeleton
$T^2$.
Normal and almost normal simple closed curves bounding
pieces of elements of
$\Sigma$
and
$\Lambda$
will partition the 2-skeleton
$T^2$
into piecewise linear regions and thus
induce a polyhedral structure on it. 
%(see figure~\ref{fig:2-skel}).
%\begin{figure}[!hbt]
%  \begin{center}
%     \vspace{2cm}
%     \Huge{CURVES IN THE BOUNDARY}
%%    \epsfig{file=protpri2.eps}
%        \caption{{\small Boundary of each tetrahedron in $T$ contains
%                         normal and almost normal simple closed curves.}}
%        \label{fig:2-skel}
%   \end{center}
%\end{figure}
We only have five non-trivial complementary regions in the boundary 
of every tetrahedron in 
$T$.
They are as in figure~\ref{fig:regions}.

So topologically we have two annuli, two twice punctured discs and
one three times punctured disc. The technical assumption on
the triangulation
$T$,
we made at the beginning of this section, implies that all these
surfaces are embedded in the 3-sphere.

We also need to take into account the discs in
$\Gamma$
which will subdivide further the polyhedral structure that 
the surfaces in figure~\ref{fig:regions}
already have. Each disc region of an element in
$\Gamma$
will give a further arc in one of the regions in 
figure~\ref{fig:regions}.
This arc will run from one normal arc in the boundary of the region
to the other. Its end points are vertices of the subdivision 
of normal and almost normal pieces of the 
2-sphere elements from
$\Lambda$.
It is worth noting that an arc in the boundary of a tetrahedron,
coming from a disc in the family
$\Gamma$,
will neither connect two segments in the 1-skeleton
$T^1$
nor will it connect a normal arc with a segment in
$T^1$.
Since a normal arc can have at most 4 vertices in its interior 
(figure~\ref{fig:superim}), it follows that we will never have to
add in more than 50 arcs per planar surface 
(adding up all the possibilities
in all the regions in figure~\ref{fig:regions}) in 
$T^2$.

We can now obtain the simplicial structure on the 2-skeleton just by
coning from one of the vertices of each disc subregion of the planar
surfaces in figure~\ref{fig:regions}. 
%But in order to make
%the shelling procedure simpler, we need to alter the polyhedral 
%structure we've constructed so far. We can achieve that by
%adding in a vertex on each
%line segment in figure~\ref{fig:regions} which is contained in
%the 1-skeleton and corresponds to 
%a disc in 
%$\Gamma$. 
%To get our simplicial structure, we cone in each 
%disc region of the 2-skeleton from one of the vertices in its
%boundary. If the disc region contains an added vertex, then 
%we make sure that this vertex is the cone point for this region.
It now follows that each surface in figure~\ref{fig:regions}
is triangulated by less than 200 2-simplices. 

%Furthermore, the 
%annular regions (the first two surfaces in the 
%figure~\ref{fig:regions}) whose line segments from
%the 1-skeleton do not correspond to any element of
%$\Gamma$,
%will be triangulated in the layered fashion. 

Now we are in the position to describe completely the subdivision
$S$
of the triangulation 
$T$
we started with.
As it was said before, the polyhedron 
$T^2\cup\Sigma\cup\Lambda\cup\Gamma$
with its simplicial structure is going to be a subcomplex 
of
$S$.
Lemma~\ref{lem:iso} tells us that the complement of
the polyhedron
$T^2\cup\Sigma\cup\Lambda\cup\Gamma$
in each tetrahedron of
$T$
is just a union of
3-balls. The boundary 2-spheres of these 3-balls 
are embedded by the assumption we made at the very beginning
of this section. They also
inherit a PL structure
from
$T^2\cup\Sigma\cup\Lambda\cup\Gamma$
which is combinatorial.
Every complementary 3-ball 
%does not contain 
%in its boundary one of the vertices 
%we had to add to the 1-skeleton
%$T^1$,
%then 
%we have the following two options.
%Either the 3-ball under consideration is a region bounded
%by two pieces belonging to two parallel elements, i.e. almost
%normal 2-spheres, in
%$\Lambda$,
%or it is not. If it is, then the intersection of its
%boundary with the tetrahedron is going to be an annulus,
%like on figure~\ref{fig:regions}, triangulated in the
%layered way. So we can just extend the layered triangulation
%over the whole 3-ball.
can thus be triangulated by adding a vertex in its interior
end coning its boundary. Since these 3-balls exhaust the
whole 3-sphere, the cones completely determine the subdivision 
$S$.
We should also note that all these coned 3-balls are in
fact shellable because their bounding 
2-spheres are triangulated in a combinatorial fashion.

%If the complementary 3-ball does contain an added vertex in its
%boundary, we don't add a new vertex in its interior.
%We just cone from the vertex in its boundary that was previously 
%added into
%the 1-skeleton. It is clear from the construction of the 
%triangulation of
%$T^2\cup\Sigma\cup\Lambda\cup\Gamma$
%that each 3-ball in the complement can contain at most one such vertex.

The rest of this section will be devoted to obtaining the subdivision
$S$
from the triangulation 
$T$
using Pachner moves. The basic tool for achieving this end will be the
procedure called {\it changing of cones}. 

\begin{figure}[!hbt]
  \begin{center}
%     \vspace{2cm}
%     \Huge{CHANGING CONES}
    \epsfig{file=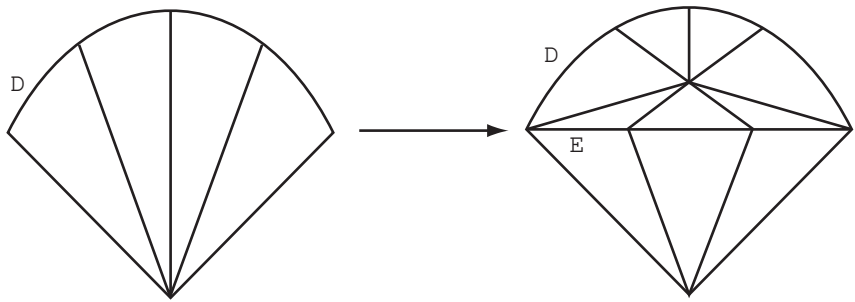}
        \caption{{\small The changing of cones.}}
        \label{fig:change}
   \end{center}
\end{figure}

Suppose we had two PL discs
$D$
and
$E$
with isomorphic simplicial structure on their boundaries. Let the union
$D\cup E$
denote the PL 2-sphere obtained by gluing the two discs together via a
simplicial isomorphism on their boundaries. What we want is an algorithm
to transform the cone on
$D$,
denoted by
$CD$,
to the union of cones 
$CE\cup C(D\cup E)$,
without changing the triangulation of
$D$.
This is described schematically by figure~\ref{fig:change}.

We have the following lemma giving a bound on the number of
Pachner moves required for changing of cones.

\begin{lem}
\label{lem:changing}
Let discs
$D$
and
$E$ 
be as above, where 
$n$
is the number of 2-simplices in
$D$
and
$m$
is the number of 2-simplices in
$E$.
Then we can perform the changing of
cones using less than
$4(n+m)$
Pachner moves. 
\end{lem}

\begin{proof}
We will divide the process into three steps. First, we glue 
a cone on the cone on the boundary of 
$D$
onto the bottom part of the boundary of
$CD$
(figure~\ref{fig:change1}).
This is a reversed process to destroying an edge which connects 
the two cone points of the bit that we glued on. It can 
therefore, by lemma~\ref{lem:crush}, be accomplished
by less than 
$(n+1)$
Pachner moves.

In the second step we perform the same move 
again, i.e. we glue the cone
$C(C(\partial D))$
onto the space we've got so far
(figure~\ref{fig:change1}). This again requires not more than 
$(n+1)$
Pachner moves.

\begin{figure}[!hbt]
  \begin{center}
%     \vspace{2cm}
%     \Huge{TWO STEPS}
    \epsfig{file=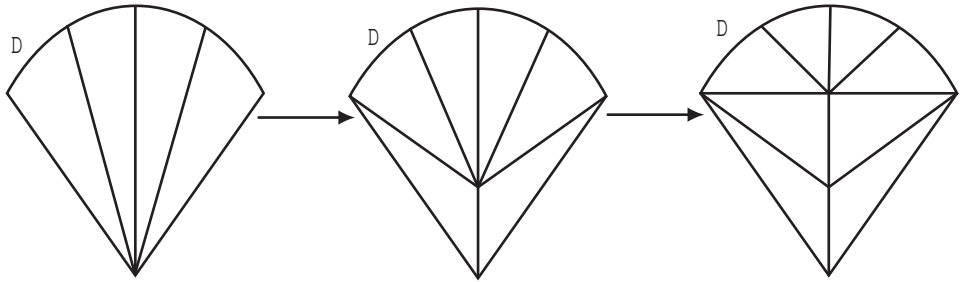}
        \caption{{\small The first two steps in the changing of cones.}}
        \label{fig:change1}
   \end{center}
\end{figure}

The space we've created can be described as a suspension of
$C(\partial D)$
glued onto the cone on the disc
$D$.
We know that we can transform the cone triangulation
of the disc
$C(\partial D)$
into the triangulation of
$E$
by using not more than 
$(n+m)$
two dimensional Pachner moves
(lemma~\ref{lem:disc}). It is also clear that in the suspension
setting, each 
$(1\p 3)$
move (or its inverse) can be realized by one
$\po$
and one
$\pt$
move. 
A
$(2\p 2)$
Pachner move can be realized by a
$\pt$
and a
$\pth$
move.
Putting all this together implies our bound.
\end{proof}

%But we also need to come up with the layered triangulations in some 
%3-ball regions of the complement of
%$T^2\cup\Sigma\cup\Lambda\cup\Gamma$.
%
%Before tackling this problem, we need to consider the following. 
%Let 
%$D$
%be a disc of the form
%$D=D'\cup S^1\times I$,
%where 
%$D'$
%is a PL 2-disc,
%$\partial D'=S^1\times0$,
%and the PL structure on
%$S^1\times I$
%is layered. 
%Again, we want to convert the PL structure of the cone
%$CD$
%into the triangulation of
%$CE\cup E\times I$,
%where 
%$E$
%is a PL disc as before and the product 
%$E\times I$
%is triangulated as follows:
%$$ E\times 0= E,\> E\times 1=D', \partial E\times I=S^1\times I$$
%and the interior of
%$E\times I$
%is layered.
%Like before, we have the following lemma.
%
%\begin{lem}
%\label{lem:changelayer}
%Let PL discs
%$D$
%and
%$E$
%be as in the lemma.
%\end{proof}

The changing of cones will help us produce all the necessary cones
in the triangulation 
$S$.
Now we have at our disposal 
all the tools required, to bound the number of Pachner moves needed
for obtaining the subdivision
$S$
from the triangulation
$T$.

The whole process will be divided into five stages. We'll start by
describing each one of them, and then we'll bound the number of
moves we made.
\begin{enumerate}
\item Add a vertex into every tetrahedron and every triangle
      of the triangulation
      $T$
      and cone.
\item Subdivide the 1-skeleton of
      $T$
      to get a subcomplex of
      $S$,
      and keep the triangulation in the 3-simplices of
      $T$
      coned.
\item Subdivide the 2-skeleton of
      $T$
      to get a subcomplex of
      $S$,
      and keep the triangulation in the 3-simplices of
      $T$
      coned.
\item Chop up tetrahedra of
      $T$
      by the appropriate normal and almost normal pieces and 
      triangulate the complementary regions by coning them
      from a point in their interior.
\item Chop up the complementary regions of 4 by length two and
      length four disc regions of elements in
      $\Gamma$. Cone the complements.       
%\item In each complementary region of the polyhedron 
%      $T^2\cup \Sigma\cup \Lambda\cup \Gamma$,
%      containing an added vertex in the 1-skeleton 
%      $T^1$,
%      destroy the edge connecting this vertex with the cone point in
%      the interior of the region.  	
\end{enumerate}

We note that step 3 can be accomplished by suspending the process in 
lemma~\ref{lem:disc}. Steps 4 and 5 are possible by 
lemma~\ref{lem:changing}.
  
Adding a vertex into each 3-simplex in
$T$
takes 
$t\po$
moves.
Adding one into a triangle of
$T$ 
takes two Pachner moves: one 
$\po$ 
move followed by a
$\pt$
move.
So step 1 amounts to 
$5t$
Pachner moves
since there are precisely 
$2t$
triangles in the triangulation
$T$.

We should note that the subdivision we get after step 1
will always satisfy the technical condition we stipulated
at the beginning of this section. This is simply
because every tetrahedron of this subdivision
contains precisely one edge from the 1-skeleton
of
$T$.
Its other edges are embedded in the 2-simplices
and in the tetrahedra of
$T$.
It is also clear that this subdivision contains
$12t$
3-simplices. So the worst case scenario would
make us do step 1 at the very beginning and then do the
simplification process (that we've been describing) on that 
subdivision.
So once we work out the bound for this simplification procedure,
we have to substitute each 
$t$
in the formula with
$12t$.

Let's go back to
the construction of the subdivision
$S$.
First we want to bound the number of vertices 
of
$S$
in each edge
of the triangulation 
$T$.
By lemma~\ref{lem:bound} it follows that there
are at most
$3\cdot 2^{300t^2}$
normal arcs in any triangle of
$T$,
coming from all elements in
$\Sigma$
and
$\Lambda$.
Since 
each normal arc contributes at most one point of intersection
with a single edge,
we will have less than
$3\cdot 2^{300t^2}$
vertices on any edge in the 1-skeleton
$T^1$.
Since there are less than 
$5t$
edges all together
(an Euler characteristic count), the number of 
vertices of the triangulation
$S$, contained in
$T^1$ 
will be bounded by
$15 t 2^{300t^2}$.

The star of any edge in
$T$
contains at most 
$2t$
3-simplices
in the subdivision we have so far. 
Creating a vertex on this edge 
can obviously be done in the following way: 
first make a
$\po$
move on one of the simplices in the star of the edge.
Then do a sequence, of
length at most
$2t-2$,
of
$\pt$
Pachner moves. 
Now the addition of the vertex can be finished off by a single
$\pth$
Pachner move. All together this procedure takes not more than
$2t$
Pachner moves.
Step 2 will thus require no more than 
$$ 30 t^2 2^{300t^2} $$
Pachner moves.

We already know that there will be at most
$3\cdot 2^{300t^2}$
normal arcs in any triangle of 
$T$.
So the number of regions in a 2-simplex in the 
2-skeleton 
$T^2$
is therefore bounded by the same number (plus one).
These regions correspond to the regions in the surfaces
from figure~\ref{fig:regions} and will 
thus be triangulated by less than 20 2-simplices. 
So any triangle in 
$T$
will be subdivided by at most
$60\cdot 2^{300t^2}$
2-simplices.
By lemma~\ref{lem:disc}, this configuration can be obtained
by 
$60\cdot 2^{300t^2}$
two dimensional Pachner moves (we should notice here
that before starting the process from lemma~\ref{lem:disc},
the triangles of 
$T$
were subdivided as cones on their boundaries).
Suspending this process and doing it
for all 
$2t$
2-simplices in
$T^2$
yields an upper bound of
$$ 3\cdot 10^2 t 2^{300t^2}$$
Pachner moves used in step 3.

The number of 3-ball regions, the elements of
$\Sigma$
and
$\Lambda$
produce in all tetrahedra of
$T$,
is equal to the number of normal and almost normal pieces
in all the 
2-spheres from
$\Sigma$
and
$\Lambda$
(plus 
$t$). So it is bounded above by
$ 3\cdot 2^{300t^2}$.
Using lemma~\ref{lem:changing}, we are going
to change the cone structure in every tetrahedron in
$T$.
This will be accomplished, step by step, starting from the
vertices of the tetrahedron and moving towards the
cone point in its interior. At each stage we have to change
a disc consisting of one of the surfaces in 
figure~\ref{fig:regions}, where all but one of its boundary 
components already have their corresponding normal 
and almost normal pieces 
glued in (that makes it a disc), to a disc
coming from the only normal or almost normal piece that 
hasn't yet been introduced. Since we want the region between
the two discs we've just described, to be coned, 
lemma~\ref{lem:changing} is precisely what is needed. 
It is also obvious that the disc
$D$
from lemma~\ref{lem:changing} will in this situation never
contain more than 800 triangles (this follows from the counts 
we did when defining the subdivision
$S$), while the disc
$E$,
which is just a normal or an almost normal piece, will
be triangulated by less than 200 2-simplices.
So in a single 3-ball region, we'll make less than 
$4\cdot(800+200)$
Pachner moves (lemma~\ref{lem:changing}).
In order to complete step 4 in all the tetrahedra
of 
$T$, we need to make
$$12\cdot 10^32^{300t^2}$$ 
Pachner moves.

The number of discs in 
$\Gamma$,
coming from a single element in 
$\Lambda$,
is bounded above by half the number of times the 2-sphere
in question intersects the 1-skeleton. We already
know that there are at most
$3\cdot 2^{300t^2}$
vertices on any edge in the 1-skeleton
of the triangulation
$T$.
Since there are less then
$5t$
edges in
$T^1$,
the number of elements in 
$\Gamma$
is bounded above by
$\frac{1}{2}15t2^{300t^2}<10t 2^{300t^2}$.

Each of the discs in
$\Gamma$
has at most
$t$
disc regions (by the assumption from the beginning of this
section),
coming from the 3-simplices in the star 
of the edge the particular disc corresponds to.
Each disc region 
is
triangulated by strictly less than 20 triangles.
A disc region in an element of 
$\Gamma$
will correspond to the disc
$E$
in lemma~\ref{lem:changing}.

The boundary of each disc region is a subcomplex in the 
boundary of a coned 3-ball. 
One of the complementary discs bounded by this
simple closed curve, in the boundary of the coned 3-ball,
will correspond to the disc 
$D$
in lemma~\ref{lem:changing}.
In the case of a disc region in an element of
$\Gamma$
having two arcs in its boundary embedded in the 
2-skeleton
$T^2$,
the disc corresponding to 
$D$
we were discussing before will contain six 2-simplices
(two in normal or almost normal pieces and four
in the 2-skeleton
$T^2$).

Let's look at the case of a disc region from an element in
$\Gamma$
that intersects the 2-skeleton of
$T$
in a single arc (i.e. the leftmost or the rightmost
region in figure~\ref{fig:disc}) and
corresponds to an elementary isotopy. The number of 
triangles of the 
complementary region (in the bounding 2-sphere) 
we are interested in will then be smaller than the sum of the
numbers of 2-simplices 
in the following surfaces: the disc in the 2-simplex of
$T$
our disc region is parallel to, 
the disc in the 2-simplex of 
$T$
containing a bounding arc of the disc region we are gluing in,
regions in at most three
normal triangles or regions in a normal triangle and a normal 
quadrilateral or regions in two normal quadrilaterals, at most
two discs contained in two distinct regions in the elements
of
$\Gamma$.
Bounds for the numbers of 2-simplices for the above surfaces are
as follows: 20, 2,
$3\cdot 30$
or
$2\cdot 30$
or
$2\cdot 30$,
$2\cdot 2$
respectively.
What happens with the disc regions 
belonging to the elements of
$\Gamma$
that come from natural isotopies?
In that case the disc 
$D$
from lemma~\ref{lem:changing} is composed of the following
surfaces: roughly a half of an almost normal octagon, three
discs contained in the 2-simplices of 
$T$,
a single normal triangle. The explicit bounds in this case are:
70, 
$3\cdot 20$,
70.

An upper bound on the sum of the numbers of triangles in 
$D$
and
$E$
will therefore always be strictly less than 
$300$
(we already know that a disc region in an element from
$\Gamma$
contains no more then 20 2-simplices).
So by lemma~\ref{lem:changing}, we can produce 
our disc region in this 3-ball by less than
$4\cdot 300$
Pachner moves. All together, we have to make less than
$$12\cdot 10^3 t^22^{300t^2}$$
Pachner moves in order to complete step 5.

%The number of vertices, not lying on any element from
%either
%$\Sigma$
%or
%$\Lambda$,
%we added to the 1-skeleton
%$T^1$,
%is equal to the number of elements in 
%$\Gamma$.
%It is therefore bounded above by 
%$100 t^3 t^{7t}$.
%Each of the 3-balls containing the edge we are about to destroy,
%is triangulated by less than 1000 tetrahedra. So,
%by lemma~\ref{lem:crush}i, each crushing move takes 
%at most a thousand
%Pachner moves. In step 6, we make less than
%$$ 10^5 t^3 2^{7t}$$
%Pachner moves.

Summing everything up,
estimating the resulting
expression and substituting
$t$
with 
$12t$
to account for
the technical assumption we made 
at the beginning of this section, we get the following proposition.

\begin{prop}
\label{prop:S}
Let 
$T$
be any triangulation of the 3-sphere and let
$t$
be the number of tetrahedra in it. Then the subdivision
$S$,
described at the beginning of this section, can be obtained from
$T$
by making less than
$ct^2 2^{dt^2}$
Pachner moves,
where the constant 
$c$
is bounded above by 
$5\cdot 10^6$
and the constant
$d$
is smaller than
$5\cdot 10^4$.
\end{prop}

\begin{center}
\section{\normalsize \scshape  CONCLUSION OF THE PROOF}
\label{sec:conc}
\end{center}

Now, we are in the position to bound the number of
Pachner moves needed to simplify any given triangulation
$T$
of the 3-sphere, down to the canonical triangulation
with only two tetrahedra. We will apply the shelling
techniques, developed in section~\ref{sec:pach}, to the
subdivision
$S$
of the triangulation
$T$,
described in section~\ref{sec:sub}.

The basic question we have to answer at this point is how many
tetrahedra do we have to shell in the simplifying process. 
Then we can estimate the number of Pachner moves needed for 
the process, using the fact that each elementary shelling
corresponds to a single Pachner move.

Let's bound first the total number of tetrahedra of
$S$.

This will be accomplished in two steps. First we count 
the number of 3-ball regions 
we coned, while constructing the
subdivision
$S$,
in all the tetrahedra of the triangulation
$T$.
The second step consists of bounding the number
of triangles in each of the boundaries of the 
3-balls mentioned above. Multiplying these two numbers 
gives our bound.

Lemma~\ref{lem:bound} implies that there are at most
$3\cdot 2^{300t^3}$
normal 
and 
almost normal 
pieces 
in all 3-simplices of
$T$,
coming from all normal and almost normal 
2-sphere in
$\Sigma\cup\Lambda$.
We know that 
each piece contains at most 200 triangles. Each planar
surface in the boundary of the tetrahedron 
(see figure~\ref{fig:regions}) contains
at most 50 arcs and is triangulated by at most 200
triangles. Each 3-ball component of the complement
of
$\Sigma\cup\Lambda$
in our tetrahedron will thus contain less than 50 disc
regions, coming from elements in
$\Gamma$.

So in all 3-simplices of
$T$
we'll have not more than
$50\cdot 3\cdot  2^{300t^2}$
3-ball regions. 
Since each disc region in any element in
$\Gamma$
contains less than 20 triangles, 1000 is surely an
upper bound on the number of triangles in the boundary
of any of the 3-ball regions.
There will therefore be at most 
$15\cdot 10^4  2^{300t^2}$
tetrahedra
in
$S$.

Combining proposition~\ref{prop:S} and the assumption
that there are precisely 
$12t$
tetrahedra in the triangulation 
$T$,
concludes the proof of the main theorem.

\vspace{5mm}
\begin{center}
\textbf{Acknowledgements}
\end{center}
\nopagebreak
I would like to thank my research supervisor Marc
Lackenby for many helpful (and enjoyable) conversations about
math and about everything else. 

\nocite{*}
\bibliographystyle{amsplain}
\bibliography{cite}

%\pagebreak

\noindent \textsc{Department of Pure Mathematics and Mathematical Statistics},
\textsc{Center for Mathematical Sciences},
\textsc{University of Cambridge}\\
\textsc{Wilberforce Road},
\textsc{Cambridge, CB3 0WB},
\textsc{UK}\\
\textit{E-mail address:} \texttt{a.mijatovic@dpmms.cam.ac.uk}
\end{document}